\documentclass[11pt]{article}
\usepackage{amsmath}
\usepackage{amssymb}
\usepackage{multirow}

{\catcode `\@=11 \global\let\AddToReset=\@addtoreset}
\AddToReset{equation}{section}

\usepackage{epsfig}
\usepackage{amsmath}
\usepackage{amssymb}
\usepackage{graphicx}
\usepackage{color}

\usepackage[latin1]{inputenc}

\newtheorem{theorem}{Theorem}[section]
\newtheorem{lemma}{\bf Lemma}[section]
\newtheorem{corollary}{\bf Corollary}[section]

\newtheorem{@definition}{\sc Definition}[section]

\newtheorem{@remark}{\sc Remark}[section]

\newtheorem{@example}{\sc Example}[section]

\newcommand{\beqn}{\begin{displaymath}}
\newcommand{\eeqn}{\end{displaymath}}
\newcommand{\beq}{\begin{equation}}  
\newcommand{\eeq}{\end{equation}}

\def\mathsf{\bf}

\def\R{\mathbb{R}}
\def\Z{\mathbb{Z}}

\def\leb{\operatorname{leb}}

\def\i{{\bf i}}
\def\d{\mathrm d}
\def\e{\mathrm e}

\def\E{\mathrm E}
\def\P{\mathrm P}

\def\text{\mbox}

\def\1{{\bf 1}}

\newcommand{\nn}{\nonumber}
\newcommand{\noi}{\noindent}

\def\limfdd{\renewcommand{\arraystretch}{0.5}
\begin{array}[t]{c}
\stackrel{\rm fdd}{\longrightarrow} \\
\end{array}\renewcommand{\arraystretch}{1}}

\def\eqfdd{\renewcommand{\arraystretch}{0.5}
\begin{array}[t]{c}
\stackrel{\rm fdd}{=} \\
\end{array}\renewcommand{\arraystretch}{1}}

\def\neqfdd{\renewcommand{\arraystretch}{0.5}
\begin{array}[t]{c}
\stackrel{\rm fdd}{\neq} \\
\end{array}\renewcommand{\arraystretch}{1}}

\newtheorem{thm}{Theorem}[section]

\newtheorem{prop}[thm]{Proposition}
\newtheorem{defn}[thm]{Definition}
\newtheorem{rem}{Remark}[section]

\begin{document}

\title{Anisotropic scaling of random grain model \\ with application
to network traffic
}
\author{Vytaut\.e Pilipauskait\.e$^{1,2}$ and Donatas Surgailis$^{1}$
\\
\small
$^1$ Vilnius University, 
Institute of Mathematics and Informatics, 08663 Vilnius, Lithuania\\
\small
$^2$ Universit\'{e} de Nantes,
Laboratoire de Math\'ematiques Jean Leray, 44322 Nantes Cedex 3, France
}
\maketitle

\begin{abstract}
We obtain a complete description of anisotropic scaling limits of random grain model on the plane with heavy tailed grain area distribution. 
The scaling limits have either independent or completely dependent increments along one or both coordinate axes and include stable, Gaussian and some `intermediate' infinitely divisible random fields.
Asymptotic form of the covariance function of the random grain model is obtained.  Application to superposed network traffic is included.
\end{abstract}

\smallskip

{\small
\noi {\it Keywords:} random grain model; anisotropic scaling; long-range dependence;
L\'evy sheet; fractional Brownian sheet; workload process
}

\vskip.7cm

\section{Introduction}

The present paper studies scaling limits of {\it random grain model}:
\begin{eqnarray}\label{rel}
X(t,s) \ := \ \sum_i \1 \big( \big( (t-x_i)/R_i^p, (s -y_i)/R_i^{1-p} \big) \in B \big),  \quad (t,s) \in \R^2,
\end{eqnarray}
where $B \subset \R^2$ (`generic grain') is a measurable bounded set of finite Lebesgue measure $\leb(B) < \infty$,
$0< p < 1$ is a shape parameter,
$\{(x_{i},y_{i}), R_i\}$ is a Poisson point process on $\R^2 \times \R_+ $ with intensity $\d x \d y F(\d r)$. We
assume that $F $ is a probability distribution on $\R_+$ having a density function $f$ such that
\begin{equation}  \label{fbeta}
f(r) \ \sim \ c_f \, r^{-1-\alpha} \quad \text{as } r \to \infty, \quad  \exists \, 1 < \alpha < 2, \ c_f >0.
\end{equation}
The sum in \eqref{rel} counts the number of uniformly  scattered and randomly dilated grains  $ (x_i,y_i) + R_i^P B$ containing $(t,s)$, where $R^P B := \{ (R^p x, R^{1-p} y) : (x,y) \in B \} \subset \R^2$ is the dilation of $B$ by factors $R^p $ and $R^{1-p}$ in the horizontal and vertical directions, respectively.  The case $p=1/2$ corresponds to uniform or isotropic dilation.
Note that the area $\leb ( R^P B ) = \leb(B) R$ of generic randomly dilated grain is proportional to $R$ and
does not depend on $p$ and has a heavy-tailed
distribution with finite mean $\E \leb ( R^P B ) < \infty $ and infinite second moment $\E \leb ( R^P B )^2 = \infty$
according to \eqref{fbeta}.
Condition \eqref{fbeta} also guarantees
that covariance of the random grain model is not integrable:
$\int_{\R^2} | {\rm Cov}(X(0,0), X(t,s)) | \d t \d s = \infty, $ see Sec.~3,
hence \eqref{rel} is a long-range dependent (LRD) random field (RF). Examples of the grain set $B$ are the unit ball and the unit square, leading respectively to the {\it random ellipses model}:
\begin{eqnarray}\label{relE}
X(t,s) \ = \ \sum_i \1 \big( (t-x_i)^2/R_i^{2p} +  (s -y_i)^2/R_i^{2(1-p)} \le 1 \big)
\end{eqnarray}
and the {\it random rectangles model}:
\begin{eqnarray}\label{relR}
X(t,s) \ = \ \sum_i \1 ( x_i < t \le x_i + R_i^p, \, y_i < s \le y_i + R_i^{1-p} ).
\end{eqnarray}
Note that the ratio of sides of a generic rectangle in \eqref{relR}
\begin{equation*}\label{RR}
\frac{R^p}{R^{1-p}}\ = \ R^{2p -1} \ \to \
\begin{cases}0, &0<p<1/2,   \\
\infty, &1/2< p < 1
\end{cases}
\qquad \text{as } R \to \infty,
\end{equation*}
implying that large rectangles are `elongated' or `flat' unless $p = 1/2$, and resulting in a strong anisotropy of \eqref{relR}. A similar observation applies to the general random grain model in \eqref{rel}.

The present paper obtains a complete description of {\it anisotropic} scaling limits
\begin{eqnarray}\label{Ylim}
a^{-1}_{\lambda,\gamma} \int_{(0,\lambda x] \times (0, \lambda^\gamma y]} (X(t,s)- \E X(t,s)) \d t \d s \ \limfdd \ V_\gamma(x,y), \qquad (x,y) \in \R^2_+, \quad \text{as } \lambda \to \infty
\end{eqnarray}
for the centered random grain model in \eqref{rel} under assumption \eqref{fbeta}, where $\{(0, \lambda x] \times (0, \lambda^\gamma y]
\subset \R^2_+, \, \lambda >0\} $ is a family of rectangles
with sides growing at possibly different rate $O(\lambda)$ and $O(\lambda^\gamma)$ and $\gamma >0$ is {\it arbitrary}. In \eqref{Ylim},  $a_{\lambda,\gamma} \to \infty$ is a normalization. See the end of this section for all unexplained notation.

\begin{center}
\begin{figure}[h]
\begin{center}
\vspace{-2.5cm}
\includegraphics[width=18cm, height=26cm]{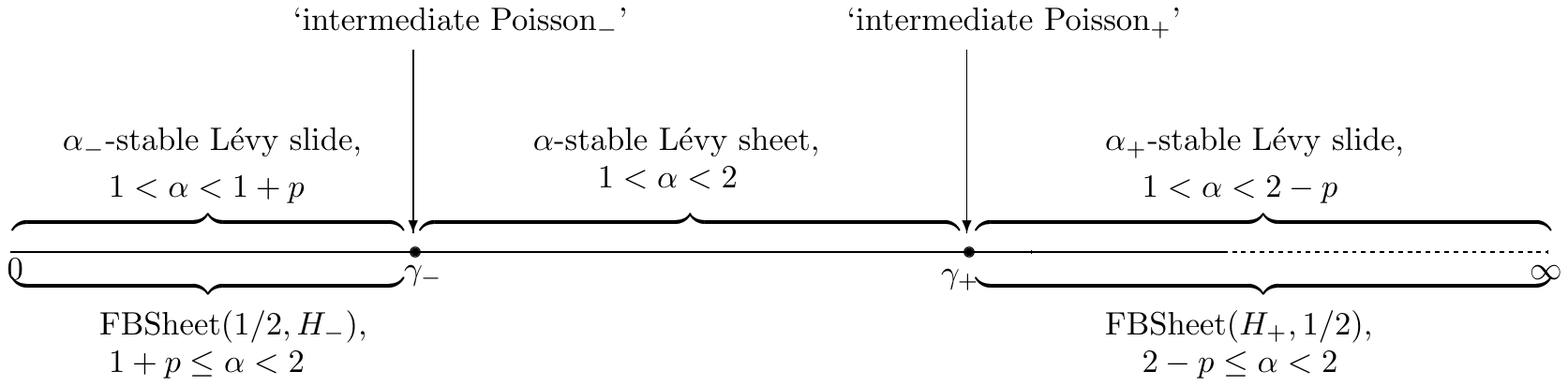}
\vspace{-18.5cm}
\caption{Scaling limits of random grain model}
\end{center}
\end{figure}
\end{center}

Our main results are summarized
in Fig.~1 which shows a  panorama of scaling limits $V_\gamma $  in \eqref{Ylim} as $\gamma $ changes between $0$ and $\infty$.
Precise formulations pertaining to Fig.~1 and the terminology therein are given in Sec.~2. Below we explain the most important facts about this diagram.
First of all note that, due to the symmetry of the random grain model in \eqref{rel}, the scaling limits in  \eqref{Ylim} are symmetric under simultaneous exchange $x \leftrightarrow y, \gamma \leftrightarrow 1/\gamma, p \leftrightarrow 1-p $ and a
reflection transformation of $B$.  This symmetry is reflected in Fig.~1, where the left region $0< \gamma \le \gamma_- $ and
the right region $\gamma_+ \le \gamma < \infty$ including
the change points of the scaling limits
\begin{equation} \label{gammapm}
\gamma_-  := \frac{1-p}{\alpha- (1-p)},    \qquad  \gamma_+  := \frac{\alpha}{p} - 1,
\end{equation}
are symmetric with respect to the above transformations.
The middle region $\gamma_- < \gamma < \gamma_+$ in Fig.~1 corresponds to an {\it  $\alpha$-stable L\'evy sheet}
defined as a stochastic integral over $(0,x]\times (0,y]$
with respect to (w.r.t.) an $\alpha$-stable random measure on $\R^2_+$.
According to Fig.~1,
for $\gamma > \gamma_+$ the scaling limits in \eqref{Ylim} exhibit a dichotomy depending on parameters $\alpha, p$, featuring a Gaussian (fractional Brownian sheet) limit for $2-p \le \alpha < 2 $, and an $\alpha_+$-stable limit for $1 < \alpha < 2-p$ with  stability parameter
\begin{equation} \label{alphaplus}
\alpha_+ := \frac{\alpha - p}{1-p} > \alpha
\end{equation}
larger than the parameter $\alpha $.
The terminology {\it $\alpha_\pm$-stable L\'evy slide} refers to a RF of the form $x L_+(y)$ or $y L_- (x)$ `sliding' linearly to zero along one of the coordinate axes, where  $L_\pm$ are $\alpha_\pm$-stable L\'evy processes (see Sec.~2 for definition).
Finally, the `intermediate Poisson' limits in Fig.~1 at $\gamma = \gamma_\pm$ are not stable although infinitely divisible RFs given by stochastic integrals w.r.t.\ Poisson random measure on $\R^2 \times \R_+$ with intensity measure $c_f \d u \d v r^{-1-\alpha} \d r $.

The results of this paper are related to works \cite{bier2010}, \cite{domb2011}, \cite{gaig2003}, \cite{kaj2007}, \cite{kaj2008}, \cite{miko2002}, \cite{pils2014}, \cite{pils2015}, \cite{ps2014}, \cite{ps2015} and others, which discuss the occurrence of different scaling regimes for various classes of LRD models, particularly, heavy-tailed duration models.
Isotropic scaling limits (case $\gamma = 1$) of random grain and random balls models in arbitrary dimension were discussed in
Kaj et al.~\cite{kaj2007}
and Bierm\'e et al.~\cite{bier2010}.
See also the monograph~\cite{lif2014} for a nice discussion of limit behavior of heavy-tailed duration models.
From an application viewpoint, probably the most interesting is the study of different scaling regimes of superposed network traffic models \cite{miko2002},  \cite{gaig2003}, \cite{kaj2007}, \cite{kaj2008}.
In these studies, it is assumed that traffic is generated by independent sources and the problem concerns the limit distribution
of the aggregated traffic as the time scale $T$ and the number of sources $M$ both tend to infinity, possibly at different rate.
The present paper extends the above-mentioned work, by considering the limit behavior
of the aggregated workload process:
\begin{eqnarray}  \label{AMKT}
A_{M,K}(Tx)&:=&\int_0^{Tx} W_{M,K}(t) \d t, \qquad \text{where} \\
W_{M, K}(t)&:=&\sum_i (R^{1-p}_i \wedge K)
\1 ( x_i < t \le x_i + R^p_i, \,
0< y_i < M ), \quad t \ge 0, \nn
\end{eqnarray}
and where $\{(x_{i},y_{i}), R_i\}$ is the same Poisson point process
as in \eqref{rel}.
The quantity $W_{M,K}(t)$ in \eqref{AMKT}
can be interpreted as the active workload at time $t$ from sources arriving at $x_i$ with $0< y_i < M$
and transmitting at rate $R^{1-p}_i \wedge K$ during
time interval $(x_i, x_i + R^p]. $  Thus, the transmission rate in \eqref{AMKT} is a (deterministic) function $(R^{p})^{(1-p)/p} \wedge K$ of the
transmission duration $R^p$ depending on parameter $0<p \le 1 $, with $0< K \le \infty $ playing the role of the maximal rate bound.
The limiting
case $p=1$ in  \eqref{AMKT} corresponds to a constant rate workload from  stationary M/G/$\infty$ queue.
Theorems~\ref{thmslow}-\ref{thminter} obtain the limit distributions of the centered and properly normalized
process  $\{A_{M,K}(Tx), \, x \ge 0\}$ with heavy-tailed distribution of $R$ in \eqref{fbeta}
when the time scale $T$, the source intensity $M$ and the maximal source rate $K$ tend jointly to infinity
so as $M = T^\gamma,  K = T^\beta $ for some $0< \gamma < \infty, 0< \beta \le \infty $. The main cases of Theorems~\ref{thmslow} and \ref{thmfast}
are summarized in Table~1.
The  workload process in  \eqref{AMKT} featuring a power-law dependence between transmission rate and duration
is closely related to the random rectangles model in  \eqref{relR}, the last fact  
being reflected in
Table~1, where most (but not all) of the limit processes can be linked to the scaling limits in Fig.~1 and where
$\gamma_+, \alpha_+$ are the same as in \eqref{gammapm}, \eqref{alphaplus}.

\begin{table}[htdp]
\begin{center}
\begin{tabular}{|l|l|l|}
\hline
\multicolumn{2}{|c|}{Parameter region} &  Limit process\\
\hline
$(1+\gamma)(1-p) <  \alpha \beta \le \infty$ & $1 < \alpha < 2$ & $\alpha$-stable L\'evy process \\
\hline
\multirow{2}{*}{$0< \alpha \beta < (1+\gamma)(1-p)$} & $1 < \alpha < 2p$ & $(\alpha/p)$-stable L\'evy process \\
\cline{2-3}
& $1 \vee 2p < \alpha < 2$ & Brownian motion\\
\hline
\end{tabular} \\
\vskip.4cm
a) Slow connection rate: $0 < \gamma < \gamma_+$
\end{center}
\end{table}

\begin{table}[htdp]
\begin{center}
\begin{tabular}{| l | l | l |}
\hline
\multicolumn{2}{|c|}{Parameter region} &  Limit process\\
\hline
\multirow{2}{*}{$0< \alpha_+ \beta < \gamma_+$} & $1 < \alpha < 2p$ & Fractional Brownian motion, $H = (3 - (\alpha/p))/2$
\\
\cline{2-3}
& $1 \vee 2p < \alpha < 2$ & Brownian motion \\
\hline
$\gamma_+ < \alpha_+ \beta < \gamma$ & \multirow{2}{*}{$1 < \alpha < 2-p$} & Gaussian line\\
\cline{1-1}\cline{3-3}
$\gamma < \alpha_+  \beta \le \infty$ & & $\alpha_+$-stable line\\
\hline
$\gamma_+ < \alpha_+ \beta \le \infty$ & $2 - p < \alpha < 2$ & Fractional Brownian motion, $H= (2-\alpha+p)/2p $ \\
\hline
\end{tabular} \\
\vskip.4cm
b) Fast connection rate: $\gamma_+ < \gamma < \infty$ \\
\caption{Limit distribution of the workload process in \eqref{AMKT} with $M = T^\gamma, K = T^\beta $
 }\label{tab:1}
\end{center}
\end{table}

The rest of the paper is organized as follows. Sec.~2 contains
rigorous formulations (Theorems~\ref{thm1}-\ref{thm5}) of the asymptotic results  pertaining to Fig.~1. Sec.~3 discusses LRD properties
and asymptotics of the covariance function of the random grain model. Sec.~4 obtains limit distributions
of the aggregated workload process in \eqref{AMKT}. 
All proofs are relegated to Sec.~5.

\smallskip

{\it Notation.} In this paper, $\limfdd$ and $\eqfdd$ denote the weak convergence and equality of finite dimensional distributions, respectively. $C$ stands for a generic positive constant which may assume different values at various locations and whose precise value has no importance.
$\R_+ := (0,\infty), \R^2_+ := (0,\infty)^2$.

\section{Scaling limits of random grain model }

We can rewrite the sum \eqref{rel} as the stochastic integral
\begin{eqnarray}\label{rel2}
X(t,s)&=&\int_{\R^2 \times \R_+}  \1 \Big( \Big( \frac{t-u}{r^p}, \frac{s -v}{r^{1-p}} \Big) \in B \Big) N (\d u, \d v, \d r), \quad (t,s) \in \R^2
\end{eqnarray}
w.r.t.\ a Poisson random measure $N(\d u, \d v, \d r)$ on $ \R^2 \times \R_+$
with intensity measure $\E N(\d u, \d v, \d r)  = \d u \d v F(\d r)$.
The integral \eqref{rel2} is well-defined and follows a Poisson distribution with mean $\E X (t,s) = \leb(B) \int_0^\infty r F (\d r)$.
The RF $X$ in
\eqref{rel2} is stationary with finite variance and the covariance function
\begin{eqnarray}\label{covX}
{\rm Cov}(X(0,0), X(t,s)) \ = \ \int_{\R^2 \times \R_+}
\1 \Big( \Big( \frac{u}{r^p}, \frac{v}{r^{1-p}} \Big) \in B,
\Big( \frac{u-t}{r^p}, \frac{v-s}{r^{1-p}} \Big) \in B \Big) \d u \d v  F(\d r).
\end{eqnarray}

Let
\begin{eqnarray}
S_{\lambda,\gamma}(x,y)&:=&\int_0^{\lambda x} \int_0^{\lambda^\gamma y} (X(t,s) - \E X(t,s)) \d t \d s  \label{Sint}\\
&=&\int_{\R^2 \times \R_+} \Big\{ \int_0^{\lambda x} \int_0^{\lambda^\gamma y} \1 \Big( \Big( \frac{t-u}{r^p}, \frac{s -v}{r^{1-p}} \Big) \in B \Big) \d t \d s  \Big\} \widetilde{N} (\d u, \d v, \d r), \quad (x,y) \in \R_+^2,\nn
\end{eqnarray}
where $\widetilde{N} (\d u, \d v, \d r) = N(\d u, \d v, \d r) - \E N(\d u, \d v, \d r)$ is the centered Poisson random measure in \eqref{rel2}.
Recall the definition of $\gamma_\pm $:
\begin{equation} \label{gammapm1}
\gamma_-  \ := \ \frac{1-p}{\alpha- (1-p)},    \qquad  \gamma_+  \ := \ \frac{\alpha}{p} - 1.
\end{equation}
The subsequent Theorems \ref{thm1}-\ref{thm5}  precise the limit RFs $V_\gamma $
and normalizations $a_{\lambda,\gamma}$ in \eqref{Ylim} for all
$\gamma >0 $ and $\alpha \in (1,2), 0<p< 1$ in Fig.~1. Throughout the paper we assume that $B$ is a bounded Borel set
whose boundary $\partial B $ has zero Lebesgue measure: $\leb (\partial B) =  0$.

\subsection{Case $\gamma_- < \gamma < \gamma_+$}

For $1 < \alpha < 2$, we introduce an $\alpha$-stable L\'evy sheet
\begin{eqnarray}\label{Lsheet}
L_\alpha (x,y) \ := \ Z_\alpha ((0,x] \times (0,y]), \quad (x,y) \in \R^2_+
\end{eqnarray}
as a stochastic integral w.r.t.\ an
$\alpha$-stable random measure  $Z_\alpha (\d u, \d v)$ on $\R^2$ with control measure $\sigma^\alpha \d u \d v$ and skewness parameter 1, where the constant $\sigma^\alpha$ is given in \eqref{const:sigma} below. Thus,
$\E \exp \{ \i \theta Z_\alpha (A) \} = \exp \{ - \leb(A) \sigma^\alpha |\theta|^\alpha (1 - \i \operatorname{sgn} (\theta) \tan(\pi \alpha/2))\}, \, \theta
\in \R$, for any Borel set $A \subset \R^2$ of finite Lebesgue measure $\leb(A) < \infty$. Note $\E Z_{\alpha}(A) = 0$.

\begin{theorem} \label{thm1}
Let $\gamma_- < \gamma < \gamma_+$, $1< \alpha < 2$. Then
\begin{eqnarray}\label{lim1}
\lambda^{-H (\gamma)}S_{\lambda,\gamma}(x,y) \ \limfdd \ L_\alpha(x,y) \quad \text{as } \lambda \to \infty,
\end{eqnarray}
where $H(\gamma) := (1+\gamma)/\alpha$  and $L_\alpha $ is an $\alpha$-stable L\'evy sheet defined in \eqref{Lsheet}.

\end{theorem}

\subsection{Cases $\gamma > \gamma_+,  1< \alpha < 2-p$ and $\gamma < \gamma_-,  1< \alpha < 1+p$ }

For $1< \alpha < 2-p$ and $1 < \alpha < 1+p$
introduce totally skewed stable L\' evy processes $\{ L_+(y), \, y \ge 0\}$ and  $\{ L_-(x), \,  x \ge 0\}$ with respective stability indices
$\alpha_\pm \in (1,2)$ defined as
\begin{equation}\label{alphapm}
\alpha_+ := \frac{\alpha -p}{1-p}, \qquad  \alpha_- := \frac{\alpha -1 + p}{p}
\end{equation}
and characteristic functions
\begin{eqnarray}\label{levyproc}
\E \exp \{\i \theta L_\pm(1) \}&:=&\exp\{-\sigma^{\alpha_\pm} |\theta|^{\alpha_\pm}(1 - \i \operatorname{sgn}(\theta) \tan(\pi \alpha_\pm/2))\},
\quad \theta \in \R,
\end{eqnarray}
where $\sigma^{\alpha_+}$ is given in \eqref{const:sigma+} and $\sigma^{\alpha_-}$ can be found by symmetry, 
see \eqref{sym} below.

\begin{theorem} \label{thm2} (i) Let $\gamma > \gamma_+$,  $1< \alpha < 2 - p$. Then
\begin{eqnarray}\label{lim2i}
\lambda^{-H(\gamma)} S_{\lambda,\gamma}(x,y) \ \limfdd \ x L_+(y) \quad \text{as } \lambda \to \infty,
\end{eqnarray}
where $H(\gamma) := 1+  \gamma/\alpha_+$  and $L_+ $ is the $\alpha_+$-stable L\' evy process defined by \eqref{levyproc}.

\smallskip

\noi (ii)  Let $0 < \gamma < \gamma_-$, $1< \alpha < 1 + p$. Then
\begin{eqnarray*}\label{lim2ii}
\lambda^{-H(\gamma)} S_{\lambda,\gamma}(x,y) \ \limfdd \ y L_-(x) \quad \text{as } \lambda  \to \infty,
\end{eqnarray*}
where $H(\gamma) := \gamma + 1/\alpha_-$  and $L_-$ is the $\alpha_-$-stable L\' evy process defined by \eqref{levyproc}.
\end{theorem}

\subsection{Cases $\gamma > \gamma_+, \ 2-p\le  \alpha < 2$ \ and \ $\gamma < \gamma_-, \ 1+p \le  \alpha < 2$ }

A (standard) fractional Brownian sheet (FBS) $B_{H_1,H_2}$ with Hurst indices $0< H_1, H_2 \le 1 $ is defined as a Gaussian process with zero mean
and covariance
\begin{eqnarray*} \label{covB}
\E B_{H_1,H_2}(x_1,y_1)B_{H_1,H_2}(x_2,y_2)&=&(1/4)(x_1^{2H_1} + x_2^{2H_1} - |x_1-x_2|^{2H_1})(y_1^{2H_2} + y_2^{2H_2} - |y_1-y_2|^{2H_2}),  \\
&&\hskip3cm (x_i,y_i) \in \R^2_+, \ i=1,2.
\end{eqnarray*}
The constants $\sigma_+$ and $\widetilde{\sigma}_+$ appearing  in Theorems \ref{thm3}~(i) and \ref{thm4}~(i) are defined in \eqref{const:sigma2_+} and \eqref{wtisigma}, respectively. The corresponding constants $\sigma_-$ and $\widetilde{\sigma}_-$ in parts (ii) of these theorems  
can be found by symmetry (see \eqref{sym}).

\begin{theorem}\label{thm3} (i) Let $\gamma > \gamma_+$, $2 - p < \alpha < 2$. Then
\begin{eqnarray}\label{lim3i}
\lambda^{-H(\gamma)} S_{\lambda,\gamma}(x,y) \ \limfdd \ \sigma_+ B_{H_+, 1/2}(x,y) \quad \text{as } \lambda \to \infty,
\end{eqnarray}
where $H(\gamma) := H_+ \,+ \,\gamma/2, \, H_+ := 1/p \, - \, \gamma_+/2 = (2-\alpha +p)/2p \in (1/2,1) $  and $ B_{H_+, 1/2}$ is
an  FBS with parameters $(H_+,1/2)$.

\smallskip

\noi (ii) Let $\gamma < \gamma_-$, $1 + p < \alpha < 2$. Then
\begin{eqnarray*}\label{lim3ii}
\lambda^{-H(\gamma)} S_{\lambda,\gamma}(x,y) \ \limfdd \ \sigma_- B_{1/2, H_-}(x,y) \quad \text{as } \lambda \to \infty,
\end{eqnarray*}
where $H(\gamma) := \gamma H_- \, + \, 1/2, \, H_- := 1/(1-p) + (1-p-\alpha)/2(1-p) \in (1/2,1)$  and $ B_{1/2,H_-}$ is
an FBS with parameters $(1/2,H_-)$.

\end{theorem}

\begin{theorem}\label{thm4} (i) Let $\gamma > \gamma_+$, $\alpha = 2 - p$. Then
\begin{eqnarray}\label{lim4i}
\lambda^{-H(\gamma)} (\log \lambda)^{-1/2} S_{\lambda,\gamma}(x,y) \ \limfdd \ \widetilde{\sigma}_+ B_{1, 1/2}(x,y) \quad \text{as } \lambda \to \infty,
\end{eqnarray}
where $H(\gamma) := 1 + \gamma/2$, $ B_{1, 1/2}$ is
an FBS with parameters $(1,1/2).$

\noi (ii) Let $\gamma < \gamma_-$, $\alpha = 1 + p$. Then
\begin{eqnarray*}\label{lim4ii}
\lambda^{-H(\gamma)} (\log \lambda)^{-1/2}  S_{\lambda,\gamma}(x,y) \ \limfdd \ \widetilde{\sigma}_- B_{1/2, 1}(x,y) \quad \text{as } \lambda \to \infty,
\end{eqnarray*}
where $H(\gamma) := \gamma  + 1/2 $ and  $ B_{1/2,1}$ is
an FBS with parameters $(1/2,1)$.
\end{theorem}

\subsection{Cases $\gamma = \gamma_\pm$  }

Define `intermediate Poisson' RFs $I_\pm = \{ I_\pm (x,y), \, (x,y) \in \R^2_+\}$ as stochastic integrals
\begin{eqnarray} \label{Inter}
I_+(x,y)&:=&\int_{\R \times (0,y] \times \R_+} \widetilde{M} (\d u, \d v, \d r) \int_{(0,x] \times \R}
\1 \Big( \Big(\frac{t-u}{r^p}, \frac{s}{
r^{1-p}} \Big) \in B \Big) \d t \, \d s,   \\
I_-(x,y)&:=&\int_{(0,x] \times \R  \times \R_+ } \widetilde{M} (\d u, \d v, \d r) \int_{\R \times (0,y]}
\1 \Big( \Big(\frac{t}{r^p}, \frac{s-v}{
r^{1-p}} \Big) \in B \Big) \d t \, \d s    \nn
\end{eqnarray}
w.r.t.\ the centered Poisson random measure $\widetilde{M} (\d u, \d v, \d r) = M (\d u, \d v, \d r) - \E M (\d u, \d v, \d r)$ on
$\R^2 \times \R_+$ with intensity measure $\E M (\d u, \d v, \d r) = c_f \d u \d v r^{-(1+\alpha)} \d r. $

\begin{prop}  \label{prop:1}
\noi (i) The RF $I_+$ in \eqref{Inter} is well-defined for
$1 < \alpha < 2$, $0<p<1$ and
$\E | I_+ (x,y) |^q < \infty$ for any $0 < q < \alpha_+ \wedge 2$.  Moreover, if $2-p < \alpha < 2$ then $\E | I_+ (x,y) |^2 < \infty$ and
\begin{equation}\label{covI}
\E I_+ (x_1, y_1) I_+ (x_2, y_2) \ = \ \sigma^2_+ \E B_{H_+, 1/2} (x_1, y_1) B_{H_+, 1/2} (x_2, y_2), \quad (x_i,y_i) \in \R^2_+, \ i = 1,2,
\end{equation}
where $\sigma_+, H_+ $ are the same as in Theorem \ref{thm3} (i).

\smallskip

\noi (ii) The RF $I_-$ in \eqref{Inter} is well-defined for $1 < \alpha < 2$, $0 < p < 1$ and
$\E | I_- (x,y) |^q < \infty$ for any $0 < q < \alpha_- \wedge 2$. Moreover,
if $1+p < \alpha < 2$ then $\E | I_- (x,y) |^2 < \infty$ and
$$
\E I_- (x_1, y_1) I_- (x_2, y_2) \ = \ \sigma_-^2 \E B_{1/2, H_-} (x_1, y_1) B_{1/2, H_-} (x_2, y_2), \quad (x_i, y_i) \in \R^2_+, \ i = 1, 2,
$$
where $\sigma_-, H_- $ are the same as in Theorem \ref{thm3} (ii).
\end{prop}

\begin{theorem}\label{thm5} (i) Let $\gamma = \gamma_+$, $1 < \alpha < 2$. Then
\begin{eqnarray}\label{lim5i}
\lambda^{-H(\gamma)} S_{\lambda,\gamma}(x,y) \ \limfdd \ I_+(x,y) \quad \text{as } \lambda \to \infty,
\end{eqnarray}
where $H(\gamma) := 1/p$ and RF $I_+$ is defined in  \eqref{Inter}.

\noi (ii) Let $\gamma = \gamma_-$, $1 < \alpha < 2$. Then
\begin{eqnarray*}\label{lim5ii}
\lambda^{-H(\gamma)}S_{\lambda,\gamma}(x,y) \ \limfdd \ I_-(x,y) \quad \text{as } \lambda \to \infty,
\end{eqnarray*}
where $H(\gamma) := \gamma_-/(1-p)$ and RF $I_-$ is defined in \eqref{Inter}.

\end{theorem}

\begin{rem} \label{rem5}
{\rm It can be easily verified that the `intermediate Poisson' RFs  $I_\pm$  in \eqref{Inter} have stationary rectangular increments
(see \cite{ps2014}, \cite{ps2015} for the definition)
and satisfy the operator self-similarity property in  \cite{bier2007}, viz.,
$\{ I_\pm (\lambda x, \lambda^{\gamma_\pm} y)\}
\eqfdd  \{\lambda^{H(\gamma_\pm)}  I_\pm (x,y) \} $  for any $\lambda >0$.  }
\end{rem}

\begin{rem} \label{rem7}
{\rm The normalizing exponent $H(\gamma) \equiv H(\gamma, \alpha, p)$ in Theorems
\ref{thm1}-\ref{thm5} is  a jointly continuous (albeit non-analytic) function
of $(\gamma, \alpha, p) \in (0,\infty) \times (1,2) \times (0,1)$.

}
\end{rem}

\begin{rem} \label{rem8}
{\rm Restriction $\alpha < 2 $ is crucial for our results. Indeed, if $\alpha > 2 $ then for
{\it any} $\gamma >0$,  $p \in (0,1)$ the
normalized integrals
$\lambda^{-(1+\gamma)/2} S_{\lambda, \gamma}(x,y)  \limfdd \sigma B_{1/2,1/2}(x,y)$  tend
to a classical Brownian sheet with variance $\sigma^2 = \leb (B)^2 \int_0^\infty r^2 F(\d r)$.
We omit the proof of the last result which follows a general scheme of the proofs
in Sec.~5.

}
\end{rem}

\section{LRD properties of
random grain model}

It is well-known that scaling limits characterize the dependence structure and large-scale properties of the underlying
random process. Anisotropic scaling of a stationary RF $Y$ on $\R^2$  as in \eqref{Ylim} with arbitrary $\gamma >0$
results in a one-dimensional family $\{V_\gamma, \gamma >0\}$ of  scaling limits and provides a more complete `large-scale summary
of $Y$' compared to the usual (isotropic) scaling with fixed $\gamma =1 $.   \cite{ps2014} observed
that for many LRD RFs $Y$ in $\Z^2$, there exists a {\it unique} point $\gamma_0 >0$ such that the scaling limits $V_\gamma \eqfdd V_\pm $ do not
depend on $\gamma $ for $\gamma < \gamma_0$ and $\gamma > \gamma_0$ and $V_+ \neqfdd V_- $. \cite{ps2014} termed
this phenomenon scaling transition (at $\gamma = \gamma_0$). The existence of scaling transition
was established for a class of aggregated nearest-neighbor autoregressive RFs \cite{ps2014} and a
natural class of Gaussian LRD RFs \cite{ps2015}. It also arises under joint temporal and contemporaneous aggregation
of independent LRD processes in telecommunication and economics, see \cite{miko2002}, \cite{gaig2003}, \cite{kaj2008},
\cite{pils2014}, \cite{pils2015}, also (\cite{ps2014}, Remark 2.3).  The results of the present work (Fig.~1)  show
a more complicated picture
with {\it two} change-points $\gamma_- < \gamma_+$ of scaling limits  which does not
fit into the definition of scaling transition in \cite{ps2014} and suggests that this concept might be more complex and
needs further studies.

One of the most common definitions of LRD property pertains to stationary random processes with non-summable (non-integrable)
autocovariance function. In the case of anisotropic RFs, the autocovariance function may decay
at different rate in different directions, motivating a more detailed classification of LRD
as in Definition \ref{lrd} below.
In this Sec.\ we also verify these LRD properties   for the random grain model in \eqref{rel}-\eqref{fbeta} and
relate them to the change
of the scaling limits or the dichotomies in Fig.~1; see Remark \ref{dich} below.

\begin{defn} \label{lrd} Let $Y = \{Y(t,s), \, (t,s) \in \R^2\}$ be a stationary RF with finite variance and nonnegative
covariance function
$\rho_Y (t,s) := {\rm Cov}(Y(0,0), Y(t,s)) \ge 0$. We say that:

\smallskip

\noi (i) $Y$ has {\it short-range dependence (SRD)}
property if $\int_{\R^2} \rho_Y(t,s) \d t \d s  < \infty $; otherwise we say that  $Y$ has {\it long-range dependence (LRD)} property;

\smallskip

\noi (ii) $Y$ has {\it vertical SRD}
property if $\int_{[-Q,Q] \times \R} \rho_Y(t,s) \d t \d s  < \infty $ for any $0<Q< \infty $; otherwise we say
that $Y$ has {\it vertical LRD} property;

\smallskip

\noi (iii) $Y$ has {\it horizontal SRD}
property if $\int_{\R \times [-Q,Q]} \rho_Y(t,s) \d t \d s  < \infty $ for any $0<Q< \infty $; otherwise we say
that $Y$ has {\it horizontal LRD} property.
\end{defn}

The main result of this Sec.\ is Theorem \ref{Xcov} providing the asymptotics of the covariance function
of the random grain model in \eqref{rel}-\eqref{fbeta} as $|t|+|s| \to \infty $
and enabling the verification of its integrability properties
in Definition \ref{lrd}. Let
$$
w := (|t|^{1/p} + |s|^{1/(1-p)})^p, \quad \text{for } (t,s) \in \R^2. 
$$
For $p=1/2$, $w$ is the Euclidean norm and $(w, {\rm arccos}(t/w))$ are the  polar coordinates  of $(t,s) \in \R^2, \, s \ge 0$.
Introduce a function $b(z), z \in [-1,1]$ by
\begin{eqnarray}\label{Leta}
b(z)
&:=&c_f \int_0^\infty \leb \Big( B \cap \Big(B + \big(z/r^{p}, (1 - |z|^{1/p})^{1-p}/ r^{1-p} \big)\Big)\Big)
r^{-\alpha} \d r,
\end{eqnarray}
playing the role of the `angular function' in the asymptotics \eqref{rhoas}. For
the random balls model  \eqref{relE} with $p=1/2 $, $b(z) $ is a constant function independent on $z$.

\begin{theorem}\label{Xcov}
Let $1< \alpha < 2$, $0< p < 1$.

\smallskip

\noi (i) The function $b(z)$ in \eqref{Leta}
is bounded, continuous and strictly positive on $[-1,1]$.

\smallskip

\noi (ii) The covariance function $\rho (t,s) := {\rm Cov}(X(0,0), X(t,s))$ in \eqref{covX} has the following asymptotics:
\begin{eqnarray}\label{rhoas}
\rho (t,s) \ \sim \ b ( \operatorname{sgn}(s) t/w ) w^{-(\alpha-1)/p} \quad \text{as } |t|+|s| \to \infty.
\end{eqnarray}
\end{theorem}

Theorem \ref{Xcov} implies the following bound for covariance function $\rho (t,s)  = {\rm Cov}(X(0,0), X(t,s))$
of the random grain model:  there exist $Q >  0$ and
strictly positive constants $0 < C_- < C_+  < \infty$ such that
 for any $|t|+|s| > Q$
\begin{equation} \label{covbdd}
C_- (|t|^{1/p} + |s|^{1/(1-p)})^{1-\alpha} \ \le \ \rho (t,s) \ \le \ C_+ (|t|^{1/p} + |s|^{1/(1-p)})^{1-\alpha}.
\end{equation}
The bounds in \eqref{covbdd} together with easy integrability properties of the function $({|t|^{1/p} + |s|^{1/(1-p)}})^{1-\alpha}$ on
$\{|t|+|s| > Q \}$ imply the following corollary.

\smallskip

\begin{corollary} The random grain model in \eqref{rel}-\eqref{fbeta} has:

\smallskip

\noi (i) LRD property for any $1 < \alpha < 2, \ 0< p < 1$;

\smallskip

\noi (ii) vertical LRD property for $1 < \alpha \le 2-p$ and vertical SRD property for
$2-p < \alpha < 2$ and any $0< p < 1$;

\smallskip

\noi (iii) horizontal LRD property for $1 < \alpha \le 1+ p$ and horizontal SRD property for
$1+ p < \alpha < 2$ and any $0< p < 1$.
\end{corollary}

\begin{rem}\label{dich}
{\rm
The above corollary indicates that the dichotomy at $\alpha = 2-p$  in Fig.~1, region $\gamma > \gamma_+$ is related to
the change from the vertical LRD  to the vertical SRD property in the random grain model.
Similarly, the horizontal transition from the LRD to the SRD explains the dichotomy at $\alpha = 1+p$ in Fig.~1, region $\gamma < \gamma_-$.
}
\end{rem}

\cite{ps2014} introduced Type I distributional LRD property for RF $Y$ with two-dimensional `time' in terms
of dependence properties of rectangular increments of scaling limits $V_\gamma, \gamma >0$.
The increment of a RF $V = \{ V(x,y), (x,y) \in \R^2_+ \}$ on rectangle $K = (u,x] \times (v,y] \subset \R^2_+$  is defined as
the double  difference  $V(K) = V(x,y) - V(u,y) - V(x,v) + V(u,v) $. Let $\ell \subset \R^2$ be a line, $(0,0) \in \ell $.
According to (\cite{ps2014}, Definition 2.2),
a RF $V = \{ V(x,y), (x,y) \in \R^2_+ \}$ is said
to have:

\begin{itemize}

\item {\it  independent rectangular increments in direction $\ell $}  if $V(K)$ and $V(K')$ are independent for any two rectangles
$K, K'  \subset \R^2_+ $ which are separated by an orthogonal line $\ell' \perp \ell $;

\item {\it invariant rectangular increments in direction $\ell $} if $V(K) = V(K') $ for any two rectangles
$K, K'$ such that $K' = (x,y) + K$ for some $(x,y) \in \ell $;

\item {\it properly dependent rectangular increments} if $V$ has neither independent nor invariant increments in
arbitrary direction $\ell $.

\end{itemize}

Further on, a stationary RF $Y$ on $\Z^2$ is said to have {\it Type I distributional LRD} (\cite{ps2014}, Definition 2.4)
if there exists a unique point $\gamma_0 >0$
such that its scaling limit $V_{\gamma_0}$ has properly dependent rectangular increments while all other scaling
limits $V_\gamma, \gamma \ne \gamma_0$  have either independent or invariant rectangular increments in some direction
$\ell = \ell (\gamma)$. The above definition trivially extends to RF  $Y$ on $\R^2$.

We end this Sec.\ with the observation that {\it all scaling limits of the random grain model in \eqref{rel}-\eqref{fbeta}  in
Theorems \ref{thm1}-\ref{thm5}
have either independent or invariant rectangular increments in direction of one or both coordinate axes. } The last
fact is immediate from stochastic integral representations in  \eqref{Lsheet},  \eqref{Inter}, the covariance function
of FBS with Hurst indices $H_1, H_2$ equal to $1 $ or $1/2$  (see also (\cite{ps2014}, Example 2.3)) and
the limit RFs in \eqref{lim2i} and \eqref{lim2ii}. We conclude that the random grain model in \eqref{rel}-\eqref{fbeta}
{\it does not have Type I distributional LRD} in contrast to Gaussian and other classes of LRD RFs discussed in
\cite{ps2014}, \cite{ps2015}. The last conclusion is not surprising since similar facts about scaling
limits of heavy-tailed duration models with one-dimensional time are well-known; see e.g. \cite{lps2005}.

\section{Limit distributions of aggregated workload process}

We rewrite the accumulated workload in \eqref{AMKT} as the integral
\begin{eqnarray}  \label{AMKT1}
A_{M,K}(Tx)&=&\int_{\R \times (0,M] \times \R_+} \Big\{(r^{1-p} \wedge K) \int_0^{Tx} \1 (u < t \le u + r^p) \d t \Big\} N(\d u, \d v, \d r),
\end{eqnarray}
where $N(\d u, \d v, \d r)$ is the same Poisson random measure on $\R^2 \times \R_+$
with intensity  $\E N(\d u, \d v, \d r) =  \d u \d  F(\d r)$ as in \eqref{rel}.
We  assume that $F(\d r)$ has a density $f(r)$ satisfying \eqref{fbeta} with $1< \alpha < 2 $ as in Sec.~2.
We let  $p\in (0,1]$ in \eqref{AMKT1} and thus the parameter may take value $p=1$ as well.
We assume that $K$ and $M $ grow with $T$ in such a way  that
\begin{equation}
M = T^\gamma, \  K = T^\beta \quad \text{for some } 0< \gamma < \infty, \ 0 < \beta \le \infty.
\end{equation}
We are interested in the limit distribution
\begin{equation}\label{Aconv}
b_T^{-1} (A_{M,K} (Tx) - \E A_{M,K} (Tx)) \ \limfdd \ {\mathcal A}(x) \quad \text{as } T \to \infty,
\end{equation}
where $b_T \equiv b_{T,\gamma,\beta} \to \infty $ is a normalization.

Recall from \eqref{gammapm} and \eqref{alphaplus} the definitions
\begin{equation*}
\gamma_+  \ = \ \frac{\alpha}{p} - 1, \qquad \alpha_+ \  = \ \frac{\alpha - p}{1-p}.
\end{equation*}
For $p=1$, let $\alpha_+  := \infty $. By assumption \eqref{fbeta}, transmission durations $R_i^p, i \in\Z $ have a heavy-tailed
distribution with tail parameter $\alpha/p  > 1$.
Following the terminology in \cite{domb2011}, \cite{gaig2003}, \cite{kaj2007}, \cite{miko2002},
the regions $\gamma < \gamma_+$, $\gamma > \gamma_+ $ and $ \gamma = \gamma_+ $
will be respectively  referred to as {\it slow connection rate},
{\it fast connection rate} and {\it intermediate connection rate}.
For each of these `regimes', Theorems \ref{thmslow},
\ref{thmfast} and \ref{thminter} detail the limit processes and normalizations
in \eqref{Aconv} depending on parameters $\beta, \alpha, p$. 
Apart from the classical Gaussian and stable processes listed in Table~1, some `intermediate' infinitely divisible processes arise. Let us introduce
\begin{eqnarray}\label{def:int_gaig}
I (x) \ := \ \int_{\R \times \R_+} \Big\{ \int_0^x 1( u < t < u + r^p ) \d t \Big\} \widetilde{\cal M} (\d u, \d r), \quad x \ge 0,
\end{eqnarray}
where $\widetilde{M} (\d u, \d r)$ is a centered Poisson random measure with intensity measure $c_f \d u r^{-(1+\alpha)} \d r$.
The process in \eqref{def:int_gaig} essentially depends on the ratio $\alpha/p$ only and is well-defined for $1 < \alpha < 2p$ and $1/2 < p \le 1$. Under the `intermediate' regime this process arises for many heavy-tailed duration models (see e.g. \cite{domb2011}, \cite{gaig2003}, \cite{kaj2008}). It was studied in detail in
\cite{gaig2006}. We introduce a `truncated' version of \eqref{def:int_gaig}:
\begin{eqnarray}\label{def:int_p}
\widehat I (x) \ := \ \int_{\R \times \R_+} \Big\{ (r^{1-p} \wedge 1) \int_0^x \1 (u < t < u + r^p) \d t \Big\} \widetilde{\cal M} (\d u, \d r), \quad x \ge 0
\end{eqnarray}
and its Gaussian counterpart
\begin{eqnarray}\label{def:int_gauss}
\widehat Z (x) \ := \ \int_{\R \times \R_+} \Big\{ (r^{1-p} \wedge 1) \int_0^x \1 (u < t < u + r^p) \d t \Big\} {\cal Z}(\d u, \d r), \quad x \ge 0,
\end{eqnarray}
where ${\cal Z} (\d u, \d r)$ is a Gaussian random measure on $\R \times \R_+$ with the same variance $c_f \d u r^{-(1+\alpha)} \d r$ as
the centered Poisson random measure  $\widetilde{\cal M} (\d u, \d r)$.  The processes in \eqref{def:int_p} and
\eqref{def:int_gauss} are well-defined for any $1< \alpha < 2$, $0< p \le 1$
and have the same covariance functions.

The RFs defined in Sec.~2 reappear in Theorems~\ref{thmslow}-\ref{thminter} for the certain grain set, namely the unit square $B=\{ (u,v) : 0 < u,v < 1 \} \subset \R^2$.
Recall that a homogeneous L\'evy process $ \{ L(x), \, x \ge 0 \} $ is completely specified by its characteristic function $\E \e^{\i \theta L(1)}, \theta \in \R $.  
A (standard) fractional Brownian motion with Hurst parameter $H \in (0,1]$ is a Gaussian process $\{B_H(x), \, x \ge 0 \}$ with zero mean and covariance
function $(1/2) (x^{2H}  + y^{2H} - |x-y|^{2H}), \, x,y \ge 0$.

\begin{theorem} \label{thmslow} {\rm (Slow connection rate.)}
 Let $0< \gamma < \gamma_+$.  The convergence in \eqref{Aconv} holds
with the limit ${\cal A}$ and
normalization $b_T = T^{\mathcal H}$ specified in (i)-(v) below.

\vskip.2cm

\noi (i) Let $(1+ \gamma) (1-p) <  \alpha \beta \le \infty$. Then $
{\mathcal H} := (1+\gamma)/\alpha $ and ${\cal A} := \{ L_\alpha (x,1), \, x \ge 0 \}$ is an $\alpha$-stable L\'{e}vy process defined by \eqref{Lsheet}.

\smallskip

\noi (ii) Let $0<  \alpha \beta < (1+ \gamma) (1-p)$ and $1 < \alpha < 2p$. Then
${\mathcal H} := \beta + (1+ \gamma) p/ \alpha $ and ${\mathcal A} := \{ L_{\alpha/p} (x), \, x \ge 0 \}$ is an $(\alpha/p)$-stable L\'{e}vy process with characteristic function given by \eqref{def:chf}.

\smallskip

\noi (iii) Let $0< \alpha \beta < (1+ \gamma) (1-p)$ and $1 \vee 2p < \alpha < 2$. Then
${\mathcal H} := (1/2)( 1 + \gamma  + \beta (2-\alpha)/(1-p) )$ and ${\mathcal A} := \{ \sigma_1 B(x), \, x \ge 0 \}$ is a Brownian motion with variance $\sigma^2_1 $ given by \eqref{const:sigma1}.

\smallskip

\noi (iv) Let $0< \alpha \beta < (1+ \gamma) (1-p)$ and $\alpha = 2p$. Then $b_T := T^{{\mathcal H}} (\log T)^{1/2}$ with
${\mathcal H} := \beta + ( 1 + \gamma)/2$ and ${\mathcal A} := \{ \widehat \sigma_1 B(x), \, x \ge 0 \}$
is a Brownian motion with variance $\widehat \sigma^2_1$ given by \eqref{const:sigma1_hat}.

\smallskip

\noi (v) Let $\alpha \beta = (1+ \gamma) (1-p)$. Then ${\mathcal H} := (1+\gamma)/\alpha$ and ${\mathcal A} := \{ \widehat L (x), \, x \ge 0 \}$ is a L\'{e}vy process with characteristic function in \eqref{def:intlevy}.
\end{theorem}

\begin{theorem} \label{thmfast} {\rm (Fast connection rate.)} Let $\gamma_+ < \gamma < \infty$.  The convergence in \eqref{Aconv} holds
with  the limit ${\cal A}$ and normalization $b_T := T^{\mathcal H}$
specified in (i)-(ix) below.

\vskip.2cm

\noi (i) Let $0 < \alpha_+\beta < \gamma_+$ and $1 < \alpha < 2p$. Then
${\cal H} := H + \beta + \gamma/2$ and ${\cal A} := \{ \sigma_2 B_H (x), \, x \ge 0 \}$ is a
fractional Brownian motion with $H = (3 - \alpha/p)/2$ and  variance $\sigma_2^2$ given by \eqref{const:sigma2}.

\smallskip

\noi (ii)  Let $0< \alpha_+ \beta < \gamma_+$ and $1 \vee 2p < \alpha < 2$. Then ${\cal H}$ and ${\cal A}$ are the same as in Theorem \ref{thmslow} (iii).

\smallskip

\noi (iii) Let $\gamma_+ < \alpha_+ \beta < \gamma$ and $1 < \alpha < 2 - p$. Then
${\cal H} := 1 + (1/2)(\gamma + \beta( 2 - \alpha - p )/(1-p))$ and ${\cal A} := \{ x Z, \, x \ge 0 \}$
is a Gaussian line  with random slope $Z \sim N(0,\sigma_3^2) $
and
$\sigma_3^2$ given in \eqref{const:sigma3}.

\smallskip

\noi (iv) Let $\gamma < \alpha_+ \beta \le \infty$ and $1 < \alpha < 2 - p$. Then
${\cal H} := 1 + \gamma/\alpha_+$ and ${\cal A} := \{ x L_+ (1), \, x \ge 0 \}$ is
an $\alpha_+$-stable line  with random slope $L_+ (1)$ having  $\alpha_+$-stable distribution defined by \eqref{levyproc}.

\smallskip

\noi (v) Let $\gamma_+ < \alpha_+ \beta \le \infty$ and $2-p < \alpha < 2$. Then
${\cal H} := H_+ + \gamma/2$ and ${\cal A} := \{ \sigma_+ B_{H_+, 1/2} (x,1), \, x \ge 0 \}$ is a fractional Brownian motion with $H = H_+ = (2 - \alpha + p)/2p$ and variance $\sigma_+^2$ given by \eqref{const:sigma2_+}.

\smallskip

\noi (vi) Let $0 < \alpha_+ \beta < \gamma_+$  and $\alpha = 2p$. Then $b_T := T^{\cal H} (\log T)^{1/2}$ with ${\cal H} := \beta + (1+ \gamma) / 2$ and
 ${\cal A} := \{ \widehat \sigma_2 B(x), \, x \ge 0 \}$ is a Brownian motion with variance  $\widehat \sigma_2^2$ in  \eqref{const:sigma_hat_2}.

\smallskip

\noi (vii) Let $\alpha_+ \beta = \gamma_+$. Then ${\cal H} := (1/2) (1 + \gamma + (2-\alpha)/p)$ and ${\cal A} := \{ \widehat Z (x), \, x \ge 0 \}$ in an intermediate Gaussian process defined by \eqref{def:int_gauss}.

\smallskip

\noi (viii) Let $\alpha_+ \beta = \gamma$ and $1 < \alpha < 2 - p$. Then ${\cal H} = 1 + \beta$ and ${\cal A} := \{ x \widehat Z, \, x \ge 0 \}$, where a slope $\widehat Z$ is a r.v. defined by \eqref{def:rv_hatZ}.

\smallskip

\noi (ix) If $\gamma_+ < \alpha_+ \beta \le \infty$ and $\alpha = 2 - p$. Then $b_T := T^{\cal H} (\log T)^{1/2}$, ${\cal H} := 1 + \gamma/2$ and ${\cal A} := \{ \widetilde \sigma_+ B_{1,1/2} (x,1), \, x \ge 0 \} = \{ x \widetilde Z, \, x \ge 0 \}$ is a Gaussian line with random slope $\widetilde Z \sim N (0,\widetilde  \sigma_+^2)$ and $\widetilde \sigma_+^2$ given by \eqref{wtisigma}.

\end{theorem}

\begin{theorem} \label{thminter}  {\rm (Intermediate connection rate.)} Let $\gamma = \gamma_+$.  The convergence in \eqref{Aconv} holds
with  the limit ${\cal A}$ and normalization $b_T := T^{\mathcal H}$
specified in (i)-(v) below.

\vskip.2cm

\noi (i) Let $0 < \alpha_+ \beta < \gamma_+$ and $1 < \alpha < 2p$. Then ${\cal H} := 1 +\beta$ and ${\cal A} := \{ I (x), \, x \ge 0 \}$ is an intermediate process defined by \eqref{def:int_gaig}.

\smallskip

\noi (ii)  Let $0 < \alpha_+ \beta < \gamma_+$ and $1 \vee 2p < \alpha < 2$. Then ${\cal H}$ and ${\cal A}$ are the same as in Theorem \ref{thmslow} (iii).

\smallskip

\noi (iii) Let $0 < \alpha_+ \beta < \gamma_+$ and $\alpha = 2p$. Then ${\cal H}$ and ${\cal A}$ are the same as in Theorem \ref{thmslow} (iv).

\smallskip

\noi (iv)  Let $\alpha_+ \beta = \gamma_+$. Then ${\cal H} := 1/p$ and ${\cal A} := \{ \widehat I (x), \, x \ge 0 \}$ is an intermediate process defined by \eqref{def:int_p}.

\smallskip

\noi (v) Let $\gamma_+ < \alpha_+ \beta \le \infty$. Then ${\cal H} := 1/p$ and ${\cal A} := \{ I_+ (x,1), \, x \ge 0 \}$ is an intermediate process defined by \eqref{Inter}.

\end{theorem}

\begin{rem} \label{rem42}
{\rm For $\gamma = \gamma_+$ we have $(1+\gamma)(1-p)/\alpha = \gamma_+/\alpha_+ = (1-p)/p$. Note that $p=1$ implies $\gamma_+ = \alpha - 1$. In this case,  Theorem \ref{thmslow} reduces to the $\alpha$-stable limit
in (i), whereas Theorem \ref{thmfast} reduces to the fractional Brownian motion limit in (v) discussed in
\cite{miko2002} and other papers.
A similar dichotomy
appears for $\beta $ close to zero and $1 < \alpha < 2p $ with the difference that $\alpha $ is now replaced
by $\alpha/p$. Intuitively, it can be explained as follows.  For small $\beta >0$, the workload process  $W_{M,K}(t)$ in
\eqref{AMKT} behaves like a constant rate process
$ K \sum_i \1(x_i <  t \le x_i + R^p_i, \, 0 < y_i < M)$   with transmission lengths $R^p_i$ that are i.i.d.\
and follow the same distribution $\P (R^p_i > r) = \P(R_i > r^{1/p}) \sim (c_f/\alpha) r^{-(\alpha/p)}, r \to \infty $ with tail
parameter $1< \alpha/p < 2$. Therefore, for small $\beta $  our results agree with \cite{miko2002},
including the Gaussian limit in Theorems \ref{thmslow} (iii) and \ref{thmfast} (ii)  arising when the $R^p_i$'s have finite variance.
}
\end{rem}

\begin{rem} \label{rem43}
{\rm As it follows from the proof,
the random line limits in Theorem~\ref{thmfast}~(iv) and~(iii) are caused by extremely long sessions starting in the past at times
$ x_i < 0$ and lasting $R^p_i = O(T^{\gamma/\gamma_+}), \, \gamma_+ < \gamma   <  \alpha_+ \beta$ or $R^p_i = O(T^{\alpha_+ \beta/\gamma_+}), \,
\gamma_+ < \alpha_+ \beta < \gamma$, respectively,
so that typically these sessions end at times $x_i + R^p_i \gg T$.

}
\end{rem}

\section{Proofs}

\subsection{Proofs of Sections 2 and 3}

Let
\begin{eqnarray*}\label{rel3}
X^*(t,s) \ := \ \int_{\R^2 \times \R_+}  \1 \Big( \Big( \frac{t-u}{r^{1-p}}, \frac{s -v}{r^{p}} \Big) \in B^* \Big) N (\d u, \d v, \d r), \quad (t,s) \in \R^2,
\end{eqnarray*}
be a `reflected' version of  \eqref{rel2}, with  $B $ replaced by $B^* := \{ (u,v) \in \R^2: (v,u) \in B \},$  $p$ replaced by
$1-p$ and the same Poisson random measure $N(\d u, \d v, \d r)$  as in \eqref{rel2}. Let $S^*_{\lambda_*,\gamma_*}(x,y) :=
\int_0^{\lambda_* x} \int_0^{\lambda_*^{\gamma_*} y} (X^*(t,s) - \E X^*(t,s)) \d t \d s, (x,y) \in \R^2_+ $ be the corresponding
partial integral in \eqref{Sint}. If $\lambda_*, \gamma_* $ are related to
$\lambda, \gamma $ as $\lambda_* =  \lambda^\gamma,  \gamma_* = 1/\gamma $ then
\begin{equation}\label{sym}
S^*_{\lambda_*,\gamma_*}(y,x) \ \eqfdd \ S_{\lambda, \gamma}(x,y)
\end{equation}
holds by symmetry property of the Poisson random measure. As noted in the Introduction,
relation \eqref{sym} allows to reduce the limits of $S_{\lambda,\gamma}(x,y)$ as $\lambda \to \infty $
and $\gamma \le \gamma_-$ to the limits of $ S^*_{\lambda_*,\gamma_*}(y,x) $ as $\lambda_* \to  \infty $
and $\gamma_* \ge \gamma_{*+}  := \alpha/(1-p) -1 $. As a consequence,
the proofs of parts (ii) of Theorems \ref{thm2}-\ref{thm5} can be omitted
since they can be deduced from parts (i)
of the corresponding statements.

The  convergence of normalized partial integrals in \eqref{Ylim}
is equivalent to the convergence of characteristic functions:
\begin{eqnarray} \label{Mchf}
\E \exp \Big\{ \i a_{\lambda,\gamma}^{-1} \sum_{i=1}^m \theta_i S_{\lambda,\gamma} (x_i,y_i) \Big\} &\to&
\E \exp \Big\{ \i\sum_{i=1}^m \theta_i V_\gamma(x_i,y_i) \Big\} \quad \text{as } \lambda \to \infty,
\end{eqnarray}
for all  $m = 1, 2, \dots$, $(x_i,y_i ) \in \R^2_+$, $\theta_i \in \R$, 
$i=1,\dots, m$. We restrict the proof of \eqref{Mchf} to one-dimensional
convergence for $m=1$, $(x,y) \in \R_+^2$ only. The general  case of  \eqref{Mchf} follows analogously.
We have
\begin{eqnarray}
W_{\lambda,\gamma}(\theta) &:=& \log \E \exp \{\i \theta a^{-1}_{\lambda,\gamma}  S_{\lambda,\gamma}(x,y)\} \nn \\
&=&\int_{\R^2 \times \R_+}
\Psi \Big( \frac{\theta}{a_{\lambda,\gamma}} \int_0^{\lambda x} \int_0^{\lambda^\gamma y}
\1 \Big( \Big(\frac{t-u}{r^p}, \frac{s -v}{r^{1-p}} \Big) \in B \Big)  \d t \d s \Big) \d u \d v f(r) \d r.\label{TH}
\end{eqnarray}
where $\Psi(z) := \e^{\i z} - 1 - {\i z}$, $z \in \R$. We shall use the following inequality:
\begin{equation}\label{ineq:psi}
| \Psi(z) | \ \le \ \min (2 |z|, z^2/2), \quad z \in \R.
\end{equation}

\bigskip

\noi {\it Proof of Theorem \ref{thm1}.} In the integrals on the r.h.s. of \eqref{TH}
we change the variables:
\begin{equation*}\label{change1}
\frac{t-u}{r^p} \ \to \ t, \quad \frac{s-v}{r^{1-p}} \ \to \ s, \quad u \ \to \ \lambda u, \quad v \ \to \ \lambda^\gamma v, \quad r \ \to \ \lambda^{H(\gamma)} r.
\end{equation*}
This yields $W_{\lambda,\gamma}(\theta)  =  \int_0^\infty g_\lambda (r) f_\lambda (r) \d r $, where
\begin{equation*}\label{flim}
f_\lambda(r) \ := \ \lambda^{(1+ \alpha)H(\gamma)} f (\lambda^{H(\gamma)} r ) \ \to \  c_f \, r^{-(1+\alpha)}, \qquad  \lambda \to \infty
\end{equation*}
according to \eqref{fbeta}, and
\begin{eqnarray*}
g_\lambda (r)&:=&\int_{\R^2}
\Psi (\theta h_{\lambda} (u,v,r) ) \d u \d v, \\
h_{\lambda} (u,v,r)
&:=&r \int_B \1(0 < u + \lambda^{-\delta_1}r^{p}t \le x, \, 0 < v + \lambda^{-\delta_2} r^{1-p} s \le y)  \d t \d s,
\end{eqnarray*}
where the exponents $\delta_1 := 1- H(\gamma) p = (\gamma_+ - \gamma)/(1+ \gamma_+) >0$, $\delta_2 := \gamma - H(\gamma)(1- p) = (\gamma - \gamma_-)/(1+ \gamma_-) >0$.
Clearly,
\begin{equation*}\label{hlim}
h_{\lambda} (u,v,r) \ \to \ \leb(B) r \, \1 (0 < u  \le x,  0 < v  \le y), \qquad  \lambda \to \infty
\end{equation*}
for any fixed $(u, v,r) \in \R^2 \times \R_+$,  $u \not\in \{0, x\}, v \not\in \{0, y\}$, implying
\begin{equation*}
g_\lambda (r) \  \to \  x y \Psi ( \theta \leb(B) r )
\end{equation*}
for any $r > 0$. Since $\int_{\R^2} h_\lambda (u,v,r) \d u \d v =  x y  r  \leb(B)$ and $h_\lambda (u,v,r) \le Cr$, the dominating
bound $| g_\lambda (r) | \le C \min (r, r^2)$ follows by \eqref{ineq:psi}.
Whence and from Lemma \ref{lemma} we conclude that
\begin{eqnarray*}
W_{\lambda,\gamma}(\theta)
&\to& W_\gamma(\theta) \ := \ xy \, c_f\int_0^\infty (\e^{\i \theta \leb(B) r} - 1 - \i \theta  \leb(B) r) r^{-(1+\alpha)} \d r. 
\end{eqnarray*}
It remains to verify that
\begin{eqnarray*}
W_\gamma(\theta) \ = \ - x y \sigma^\alpha |\theta|^\alpha \big( 1 - \i \operatorname{sgn} (\theta) \tan(\pi \alpha/2) \big)  \ = \
\log \E \exp \{\i \theta L_\alpha (x,y) \},
\end{eqnarray*}
where
\begin{equation}\label{const:sigma}
\sigma^\alpha \ := \ c_f \leb(B)^\alpha  \cos (\pi \alpha/2) \Gamma (2-\alpha)/ \alpha (1- \alpha).
\end{equation}
This proves the one-dimensional convergence in \eqref{lim1} and Theorem  \ref{thm1}, too.
\hfill $\Box$

\bigskip

\noi {\it Proof of Theorem \ref{thm2}.} In \eqref{TH}, change the variables as follows:  
\begin{equation}\label{change2}
t  \ \to  \ \lambda t, \quad s-v \ \to \ \lambda^{(1-p)\gamma/(\alpha -p)} s,
\quad u \ \to \ \lambda^{p \gamma/(\alpha -p)}u, \quad v \ \to \ \lambda^\gamma v, \quad r \ \to \ \lambda^{\gamma/(\alpha -p)} r. 
\end{equation}
This yields $W_{\lambda,\gamma}(\theta)  =  \int_0^\infty g_\lambda (r) f_\lambda (r) \d r $, where
\begin{equation}\label{flim2}
f_\lambda(r) \ := \ \lambda^{(1+\alpha)\gamma/(\alpha -p)} f ( \lambda^{ \gamma/(\alpha -p) }r ) \ \to \  c_f \, r^{-(1+\alpha)}, \qquad  \lambda \to \infty
\end{equation}
and $g_\lambda (r):= \int_{\R^2}
\Psi (\theta h_{\lambda} (u,v,r) ) \d u \d v $ with
\begin{eqnarray} \label{h2}
h_{\lambda} (u,v,r)  
&:=&\int_0^x \d t \int_{\R} \1\Big( \Big(\frac{\lambda^{-\delta_1} t - u}{r^{p}}, \frac{s}{r^{1-p}}\Big) \in B \Big)
\1( 0 < v +   \lambda^{-\delta_2} s < y) \d s,
\end{eqnarray}
where
$\delta_1 := p\gamma/(\alpha -p) -1 = (\gamma - \gamma_+)/\gamma_+ >0$,
$\delta_2 := \gamma (\alpha -1)/(\alpha -p) >0$.
Let
$B(u) := \{ v \in \R: (u,v) \in B\}$ and write $\leb_1(A)$ for the Lebesgue measure
of a set $A \subset \R$. By the dominated convergence theorem,
\begin{eqnarray}
h_\lambda (u,v,r)\ \to\ 
h(u,v,r)&:=&
x\,\1(0 < v < y)  \int_{\R} \1\Big( \Big(\frac{-u}{r^p}, \frac{s}{r^{1-p}} \Big) \in B \Big)\d s \label{h2lim}\\
&=&x \,\1(0 < v < y) r^{1-p} \leb_1 ( B(-u/r^{p}) )\nn
\end{eqnarray}
for any $(u,v,r) \in \R^2\times \R_+, v \not\in \{0, y\}$, implying
\begin{eqnarray*}
g_\lambda (r)&\to&g(r) \ := \
\int_{\R^2} \Psi (\theta h (u,v,r) ) \d u \d v \ = \ y \,r^{p}
\int_{\R} \Psi\big(\theta x r^{1-p}  \leb_1(B(u))\big) \d u
\end{eqnarray*}
for any $r>0$. Indeed, 
since $B$ is bounded, for fixed $r>0$  the function $(u,v) \mapsto h_\lambda (u,v,r)$ has a bounded support uniformly
in $\lambda \ge 1$. Therefore it is easy to verify domination criterion for the above convergence.
Combining $h_\lambda (u,v,r) \le C r^{1-p}$ with $\int_{\R^2} h_\lambda (u,v,r) \d u \d v = x y  r \leb(B)$
gives $|g_\lambda (r) | \le C \min(r, r^{2-p})$ by \eqref{ineq:psi}.
Hence and
by Lemma~\ref{lemma}, $W_{\lambda,\gamma}(\theta) \to  W_\gamma(\theta) := c_f \int_0^\infty g(r) r^{-(1+\alpha)} \d r$.
By change of variable, the last integral can be rewritten as
\begin{eqnarray*}
 W_\gamma(\theta)&=&c_f \,y \, x^{\alpha_+}(1-p)^{-1}\int_\R  \leb_1(B(u))^{\alpha_+} \d u \int_0^\infty
 (\e^{\i \theta w} - 1 - \i \theta w ) w^{-(1+\alpha_+)} \d w\\
&=&- ( y \, x^{\alpha_+} ) \sigma^{\alpha_+} |\theta|^{\alpha_+} \big(1 - \i \operatorname{sgn} (\theta) \tan(\pi \alpha_+/2) \big) \ = \ \log \E \exp \{\i \theta x L_{+} (y) \},
\end{eqnarray*}
where
\begin{eqnarray}\label{const:sigma+}
\sigma^{\alpha_+} \ := \  \frac{c_f \Gamma (2 - \alpha_+)\cos(\pi \alpha_+/2)}{(1-p) \alpha_+ (1-\alpha_+) } 
\int_\R  \leb_1(B(u))^{\alpha_+} \d u,
\end{eqnarray}
thus completing the proof of one-dimensional convergence in \eqref{lim2i}. 
Theorem~\ref{thm2} is proved.
\hfill $\Box$

\bigskip

\noi {\it Proof of Theorem \ref{thm3}.} 
In \eqref{TH}, change the variables as follows: 
\begin{equation}\label{change3}
t \ \to \  \lambda t, \quad  s - v \ \to \ \lambda^{(1/p)-1} s, \quad u \ \to \ \lambda u, \quad v \ \to \ \lambda^\gamma v, \quad r \ \to \ \lambda^{1/p} r.
\end{equation}
We get $W_{\lambda,\gamma}(\theta) = \int_0^\infty g_\lambda (r) f_\lambda (r) \d r$, where
\begin{equation} \label{fg3}
f_\lambda (r) \ := \ \lambda^{(1+\alpha)/p} f ( \lambda^{1/p} r ), \quad g_\lambda (r) \ := \ \int_{\R^2} \lambda^{2 ( H(\gamma)-1/p)} \Psi ( \theta \lambda^{(1/p)-H(\gamma)} h_\lambda (u,v,r) ) \d u \d v,
\end{equation}
with
\begin{eqnarray} \label{h3}
h_\lambda (u,v,r) &:=& \int_0^x \d t \int_{\R} \1 ( 0 < v + \lambda^{-\delta} s < y ) \1 \Big( \Big( \frac{t-u}{r^p} , \frac{s}{r^{1-p}} \Big) \in B \Big)
\d s  \nn \\
&\to&\1 ( 0 < v < y )\int_0^x \d t \int_{\R} \1 \Big( \Big( \frac{t-u}{r^p} , \frac{s}{r^{1-p}} \Big) \in B \Big) \d s \nn \\
&=&\1 (0 < v  < y)\, r^{1-p} \int_0^x \leb_1 (B((t - u)/r^{p}))\d t \ =: \ h(u,v,r), \qquad \lambda \to \infty,
\end{eqnarray}
for all $(u,v,r) \in \R^2 \times \R_+$, $v \not \in \{ 0, y \}$, since $\delta := 1 + \gamma - (1/p) > 0.$
Note that $2(H(\gamma)-1/p) = \gamma - \gamma_+ > 0$ and hence
$$
\lambda^{2 (H(\gamma) - 1/ p)} \Psi ( \theta \lambda^{(1/p) - H(\gamma)} h_\lambda (u,v,r) ) \ \to \ - (\theta^2/2) h^2(u,v,r), \qquad \lambda \to \infty.
$$
Next, by the dominated convergence theorem 
$$
g_\lambda (r) \ \to \ g(r) : = - \frac{\theta^2}{2} \int_{\R^2} h^2(u,v,r) \d u \d v
$$
for any $r>0$.
Using $\int_{\R^2} h_\lambda (u,v,r) \d u \d v = xy \leb(B) r$ and $h_\lambda (u,v,r) \le C \min ( r^{1-p}, r )$ similarly as in the proof of
Theorem \ref{thm2}
we obtain $|g_\lambda (r)| \le C \int_{\R^2} h^2_\lambda (u,v,r) \d u \d v \le C \min (r^{2-p}, r^2)$. 
Then by Lemma \ref{lemma},
$$
W_{\lambda,\gamma}(\theta) \ \to \ W_{\lambda}(\theta) \ := \ c_f \int_0^\infty g(r) r^{-(1+\alpha)} \d r \ =  \ - (\theta^2/2) \sigma^2_+ x^{2 H_+} y,
$$
where
\begin{eqnarray}\label{const:sigma2_+}
\sigma^2_+ \ := \ c_f \int_\R \d u \int_0^\infty  \Big( \int_0^1 \leb_1 (B((t - u)/r^{p})) \d t  \Big)^2  r^{1-\alpha -2p} \d r,
\end{eqnarray}
where the last integral converges. (Indeed, since  $u \mapsto \leb_1 (B(u)) = \int 1((u,v) \in B) \d v $ is a bounded function with compact support, the inner integral in \eqref{const:sigma2_+} does not exceed $C (1 \wedge r^p) \1 (|u| < K(1+ r^p) )$ for some $C, K >0$ implying
$\sigma^2_+ \le C \int_0^\infty (1 \wedge r^p)^2 (1 + r^p) r^{1-\alpha -2p} \d r < \infty $ since $2-p < \alpha < 2 $.)
This ends the proof of one-dimensional convergence in \eqref{lim3i}. 
Theorem~\ref{thm3} is proved.
\hfill $\Box$

\bigskip

\noi {\it Proof of Theorem \ref{thm4}.} After  the same change of variables as in \eqref{change2}, viz.,
\begin{equation*}\label{change4}
t  \ \to  \ \lambda t, \quad s-v \ \to \ \lambda^{\gamma/2} s,
\quad u \ \to \ \lambda^{p \gamma/2(1 -p)}u, \quad v \ \to \ \lambda^\gamma v, \quad r \ \to \ \lambda^{\gamma/2(1 -p)} r,
\end{equation*}
we obtain $W_{\lambda,\gamma}(\theta) = \int_0^\infty g_\lambda (r) f_\lambda (r) \d r $ with $f_\lambda (r)$ as in  \eqref{flim2}
and $g_\lambda (r) := \int_{\R^2} \Psi (\theta (\log \lambda)^{-1/2} h_\lambda (u,v,r)) \d u \d v $,
where
\begin{eqnarray*} \label{h4}
h_{\lambda} (u,v,r)
&:=&\int_0^x \d t \int_{\R} \1\Big( \Big(\frac{\lambda^{-\delta_1} t - u}{r^{p}}, \frac{s}{r^{1-p}}\Big) \in B \Big)
\1( 0 < v +   \lambda^{-\delta_2} s < y) \d s, 
\end{eqnarray*}
$\delta_1 := p\gamma/2(1-p) - 1 = (\gamma - \gamma_+)/\gamma_+ >0$, $\delta_2 := \gamma/2 >0 $ are the same 
as in \eqref{h2}  and 
\begin{eqnarray*}
h_\lambda (u,v,r) \ \to \ h(u,v,r)  \ := \
x \, \1(0 < v < y)  \int_{\R} \1\Big( \Big(\frac{-u}{r^p}, \frac{s}{r^{1-p}} \Big) \in B \Big)\d s
\end{eqnarray*}
c.f. \eqref{h2lim}. Below we prove that the main contribution to the limit of $W_{\lambda,\gamma}(\theta)$ comes from
the interval $\lambda^{-\delta_1/p} < r < 1 $, namely, that $W_{\lambda,\gamma}(\theta) - W^0_{\lambda,\gamma}(\theta) \to 0$,
where
\begin{eqnarray}\label{W0}
W^0_{\lambda,\gamma}(\theta)
\ :=\ \int_{\lambda^{-\delta_1/p}}^{1} g_\lambda (r) f_\lambda (r) \d r
&\sim&-\frac{\theta^2}{2} \, \frac{c_f}{\log \lambda}  \int_{\lambda^{-\delta_1/p}}^1
\frac{\d r}{r^{3-p}}  \int_{\R^2} h^2(u,v,r) \d u \d v  \\
&=&- \frac{\theta^2}{2} x^2 y c_f\int_\R (\leb_1 (B(u)))^2 \d u \, \frac{1}{\log \lambda}  \int_{\lambda^{-\delta_1/p}}^1
r^{-1} \d r \nn \\
&=&- \frac{\theta^2}{2}\widetilde \sigma^2_+  x^2 y  \ =: \  W_\gamma(\theta), \nn
\end{eqnarray}
where
\begin{equation}\label{wtisigma}
\widetilde \sigma^2_+ \ := \ \frac{c_f (\gamma-\gamma_+)} {2(1-p)} \int_\R \leb (B \cap (B + (0,u))) \d u 
\end{equation}
and where we used the fact that 
$\int_{\R^2} h^2(u,v,r) \d u \d v = x^2y r^{2-p} \int_\R \leb_1 (B(u))^2 \d u =
x^2y r^{2-p}  \int_\R \leb (B \cap (B + (0,u))) \d u $.

Accordingly, write $W_{\lambda,\gamma}(\theta) =  W^0_{\lambda,\gamma}(\theta) + W^-_{\lambda,\gamma}(\theta) + W^+_{\lambda,\gamma}(\theta)$,
where $W^-_{\lambda,\gamma}(\theta) := \int_0^{\lambda^{-\delta_1/p}} g_\lambda (r) f_\lambda (r) \d r $ and
$W^+_{\lambda,\gamma}(\theta) := \int_1^{\infty} g_\lambda (r) f_\lambda (r) \d r$ are remainder terms. Indeed,
using \eqref{ineq:psi} and 
\begin{equation}\label{hin}
\int_{\R^2} h_\lambda (u,v,r) \d u \d v \ = \ x y r \leb (B), \qquad  h_\lambda (u,v,r) \ \le \ C (\lambda^{\delta_1} r)\wedge r^{1-p}.
\end{equation}
it follows that
\begin{eqnarray*}
|W^+_{\lambda,\gamma}(\theta)|&\le&\frac{C}{ (\log \lambda)^{1/2}} \int_{1}^\infty \frac{\d r}{r^{3-p}} \int_{\R^2} h_\lambda(u,v,r) \d u \d v
 \ = \  O( (\log \lambda)^{-1/2}) \ = \  o(1).  \label{Wplus}
\end{eqnarray*}
Similarly,
\begin{eqnarray*}
|W^-_{\lambda,\gamma}(\theta)|&\le&\frac{C\lambda^{\delta_1}}{\log \lambda} \int_0^{\lambda^{-\delta_1/p}}  r f_\lambda(r)  \d r  \int_{\R^2} h_\lambda(u,v,r) \d u \d v \
\le\ \frac{C \lambda^{\delta_1} }{\log \lambda}
\int_0^{\lambda^{-\delta_1/p}}  r^2 f_\lambda(r)  \d r \nn \\
&=&\frac{C}{ \lambda  \log \lambda}
\int_0^{\lambda^{1/p}}  r^2 f(r)  \d r \ = \ O((\log  \lambda)^{-1}) \ = \  o(1).  \label{Wmin}
\end{eqnarray*}
since $\delta_1 = p\gamma/2(1-p) -1$. 

Consider the main term $W^0_{\lambda,\gamma}(\theta)$ in \eqref{W0}. Let
$\widetilde W_{\lambda,\gamma}(\theta) := -\frac{\theta^2}{2\log \lambda}  \int_{\lambda^{-\delta_1/p}}^{1}
f_\lambda (r)
\d r  \int_{\R^2} h^2_\lambda(u,v,r) \d u \d v $. Then using \eqref{hin} and
$|\Psi (z) + z^2/2| \le |z|^3/6 $ we obtain
\begin{eqnarray*}
|W^0_{\lambda,\gamma}(\theta)- \widetilde W_{\lambda,\gamma}(\theta)|
&\le&\frac{C}{(\log \lambda)^{3/2}}
\int_{\lambda^{-\delta_1/p}}^1  r^{2-2p} f_\lambda(r)  \d r  \int_{\R^2} h_\lambda(u,v,r) \d u \d v \\
&\le&\frac{C}{(\log \lambda)^{3/2}}
\int_{\lambda^{-\delta_1/p}}^1  r^{3-2p} f_\lambda(r)  \d r  \\
&\le&\frac{C}{ (\log \lambda)^{3/2}}
\int_{0}^1  r^{-p}  \d r  \ = \  O((\log \lambda)^{-3/2}) \ = \  o(1).
\end{eqnarray*}
Finally, it remains to estimate the difference $|\widetilde W_{\lambda,\gamma}(\theta)  - W_\gamma(\theta)|
\le C(J'_\lambda + J''_\lambda)$, where 
\begin{eqnarray*}
J'_\lambda&:=&\frac{1}{\log \lambda}  \int_{\lambda^{-\delta_1/p}}^1
f_\lambda (r) \d r  \int_{\R^2} |h^2_\lambda(u,v,r) -  h^2(u,v,r)|  \d u \d v, \\
J''_\lambda&:=&\frac{1}{\log \lambda}  \int_{\lambda^{-\delta_1/p}}^1
r^{2-p} |f_\lambda (r) - c_f r^{p-3} | \d r.
\end{eqnarray*}
Let
\begin{eqnarray*} \label{h5}
\widetilde h_{\lambda} (u,v,r)
&:=&x \int_{\R} \1\Big( \Big(\frac{- u}{r^{p}}, \frac{s}{r^{1-p}}\Big) \in B \Big)
\1( 0 < v +   \lambda^{-\delta_2} s < y) \d s.
\end{eqnarray*}
Then $J'_\lambda \le  J'_{\lambda 1} + J'_{\lambda 2}$, where
$J'_{\lambda 1} :=  (\log \lambda)^{-1}  \int_{\lambda^{-\delta_1/p}}^{1}
f_\lambda (r) \d r  \int_{\R^2} |h^2_\lambda(u,v,r) -  \widetilde h^2_\lambda (u,v,r)|  \d u \d v $ and 
$J'_{\lambda 2} := $   $(\log \lambda)^{-1}  \int_{\lambda^{-\delta_1/p}}^{1}
f_\lambda (r) \d r  \int_{\R^2} |\widetilde h^2_\lambda(u,v,r) -  h^2(u,v,r)|  \d u \d v$.
Using the fact that $B$ is a bounded set with $\leb (\partial B) = 0 $  we get that
\begin{eqnarray*}
\int_{\R^2} |h_\lambda(u,v,r) - \widetilde h_\lambda (u,v,r)| \d u \d v
&\le&yr \int_0^x \d t \int_{\R^2} \Big|\1\Big( \Big( \frac{\lambda^{-\delta_1} t}{r^p} - u, s \Big) \in B \Big)
-  \1\big( (-u,s) \in B \big)\Big| \d u \d s \\
&\le&r \epsilon (1/\lambda^{\delta_1}r^{p}),
\end{eqnarray*}
where $\epsilon(z), z \ge 0$ is a bounded function with $\lim_{z \to 0} \epsilon (z) = 0$.  
We also  have
$ h_\lambda(u,v,r) + \widetilde h_\lambda (u,v,r) \le C
r^{1-p} $ as in \eqref{hin}. 
Using these bounds together with  $f_\lambda (r) \le C r^{p-3}, r > \lambda^{-\delta_1/p}$
we obtain
\begin{eqnarray*}
J'_{\lambda 1} \log \lambda 
\ \le \ C \int_{\lambda^{-\delta_1/p}}^{1}
\epsilon (1/\lambda^{\delta_1}r^{p}) r^{-1} \d r
\ = \ C \int_{\lambda^{-\delta_1}}^1  \epsilon (z) z^{-1} \d z \ = \  o(\log \lambda), 
\end{eqnarray*}
proving  $J'_{\lambda 1} \to 0$ as $\lambda \to \infty $.  In a similar way, using
$\int_{\R^2} |\widetilde h_\lambda(u,v,r) - h(u,v,r)| \d u \d v
\le xr \int_{\R^3} \1((- u, s) \in B) | \1(0 < v + \lambda^{-\delta_2} r^{1-p} s < y) 
-  \1(0< v < y) | \d u \d v \d s \le C r^{2-p}\lambda^{-\delta_2}    $ we obtain
$J'_{\lambda 2} \log\lambda  \le C \lambda^{-\delta_2} \int_0^1 r^{-p} \d r
=  O(\lambda^{-\delta_2})$, proving  $J'_{\lambda 2} \to 0$  and hence $J'_{\lambda} \to 0$. 
Finally,  $J''_{\lambda} = (\log \lambda)^{-1} \int_{\lambda^{1/p}}^\infty r^{2-p} |f(r) - c_f r^{p-3}| \d r \to 0$
follows from  \eqref{fbeta}. This proves the limit 
$ \lim_{\lambda \to \infty} W_{\lambda,\gamma}(\theta) = W_\gamma(\theta)
= -(\theta^2/2)\widetilde \sigma^2_+  x^2 y $ for any $\theta \in \R $, or
one-dimensional convergence in \eqref{lim4i}. Theorem \ref{thm4}  is proved.
 \hfill $\Box$

\bigskip

\noi {\it Proof of Proposition \ref{prop:1}.} We use well-known properties of Poisson stochastic integrals
and inequality  (3.3) in \cite{pils2014}. Accordingly, $I_+(x,y)$ is well-defined and satisfies
$\E |I_+(x,y)|^q \le 2 J_q(x,y) \, (1 \le q \le 2)$ provided 
\begin{eqnarray*}
J_q(x,y)&:=&c_f\int_0^\infty r^{-(1+\alpha)} \d r \int_{\R \times (0,y]} \d u \d v \Big| \int_{(0,x] \times\R} \1 \Big( \Big( \frac{t-u}{r^p}, \frac{s}{r^{1-p}} \Big) \in B \Big) \d t \d s \Big|^q  \\
&=&c_f y \int_0^\infty r^{q(1-p)  -(1+\alpha)} \d r 
\int_{\R} \d u \Big|\int_0^x \leb_1 \Big( B \Big( \frac{t-u}{r^p} \Big) \Big) \d t\Big|^q \ < \ \infty.
\end{eqnarray*}
Split $J_q(x,y) = c_f y [\int_0^1 \d r + \int_1^\infty ] \cdots \d r =: c_f y [J' + J''] $. Then
$J'' \le C  \int_1^\infty r^{q(1-p)  -(1+\alpha)} \d r \int \1(|u| \le C r^p) \d u \le C \int_1^\infty r^{q(1-p)  - (1+ \alpha) +p} \d r
< \infty $ provided $q < (\alpha -p)/(1-p).$  Similarly, $J' \le C  \int_0^1 r^{q(1-p)  -(1+\alpha)} \d r |\int \1(|t| \le C r^p ) \d t |^q
\le C \int_0^1 r^{q(1-p)  - (1+ \alpha) + qp} \d r
< \infty$ provided $\alpha < q$.  Note that $\alpha < (\alpha -p)/(1-p) \le 2$ for $1< \alpha \le 2-p$ and  
$(\alpha -p)/(1-p) > 2 $ for $2-p < \alpha <2$.   
Relation \eqref{covI} follows from \eqref{covB}  and
$J_2(x,y) = \sigma^2_+ y x^{2H_+}$ by a change of variables. This proves part (i). The proof of part (ii)
is analogous. \hfill $\Box$

\bigskip

\noi {\it Proof of Theorem \ref{thm5}.} Using the change of  variables as in  \eqref{change3} we get
$W_{\lambda,\gamma}(\theta) = \int_0^\infty g_\lambda (r) f_\lambda (r) \d r$ with the same  $f_\lambda (r), g_\lambda (r)$ as
in \eqref{fg3}  and $h_\lambda (u,v,r)$ satisfying \eqref{h3}. (Note $H(\gamma) = H(\gamma_+) = 1/p$ hence 
$\lambda^{H(\gamma_+) - (1/p)} =1 $ in  the definition
of $g_\lambda (r)$ in \eqref{fg3}.) Particularly, $\Psi (\theta h_\lambda (u,v,r)) \to
\Psi (\theta h(u,v,r))$ for any $(u,v,r) \in \R^2 \times \R_+$, $v \not \in \{ 0, y \}$. Then 
$g_\lambda (r) \to g(r) := \int_{\R^2} \Psi (\theta h(u,v,r)) \d u \d v$ follows by the dominated
convergence theorem. Using $\int_{\R^2} h_\lambda (u,v,r) \d u \d v = x y r \leb (B)$ and $h_\lambda (u,v,r) \le C r$ we obtain 
$|g_\lambda (r)| \le C \min (r, r^2)$  and hence  $W_{\lambda,\gamma}(\theta) \to \int_0^\infty g(r) r^{-(1+\alpha)} \d r = 
\log \E \exp \{ \i \theta I_+ (x,y) \}$, proving the one-dimensional convergence in \eqref{lim5i}. 
The proof of Theorem~\ref{thm5} is complete.
\hfill $\Box$

\bigskip

\noi {\it Proof of Theorem \ref{Xcov}.} (i) Write $D_r(x,y) := \{ (u,v)\in \R^2: (u-x)^2 + (v-y)^2 \le r^2 \}$ for a ball
in $\R^2$ centered at $(x,y)$ and having radius $r$.
Recall that $B$ is bounded. Note that $\inf_{z\in [-1,1]} ( |z| / r^p + (1 - |z|^{1/(p-1)})^{1-p} / r^{1-p} ) \ge c_0 \min (r^{-p}, r^{-(1-p)})$ for some constant $c_0 >0$. Therefore, there exists $r_0 > 0$ such that for  all $0<r < r_0$ the intersection
$B_{z,r} :=
B \cap \big(B + \big(z/r^{p}, (1 - |z|^{1/p})^{1-p}/ r^{1-p} \big)\big) = \emptyset $ in \eqref{Leta}.
Hence $b (z) \le C < \infty$ uniformly in $z \in [-1,1]$.

Let $(x,y) \in B\setminus \partial B$. Then $D_{2r}(x,y) \subset B$ for all $r < r_0$ and some $r_0 >0$.
If we translate $B$ by distance $r_0$ at most, the translated set 
still contains the ball $D_{r_0}(x,y)$.
Since $\sup_{z \in [-1,1]} ( |z| /r^p + (1-|z|^{1/p})^{1-p} / r^{1-p} ) \le 2 \max (r^{-p}, r^{-(1-p)})$, there exists  $r_1 >0$ 
for which $ \inf_{r > r_1}  \leb (B_{z,r}) \ge \pi r_0^2, $ proving   $\inf_{z \in [-1,1]} b (z) >0$.  The continuity
of $b(z)$ follows from the above argument and the continuity of the mapping $z \mapsto \leb(B_{z,r}): [-1,1] \to \R_+$, for each
$r >0$.

\smallskip

\noi (ii) Let $s \ge 0$. 
In the integral \eqref{covX} we change the variables: $u \to r^p u$, $v \to r^{1-p} v$, $r \to w^{1/p} r$.
Then
$$
\rho (t,s) \ = \  w^{-(\alpha-1)/p} \int_0^\infty \leb (B_{t/w,r}) f_w (r) r \d r,
$$
where $f_w(r) := w^{(1+ \alpha)/p} f(w^{1/p} r) \to c_f \, r^{-(1+\alpha)}$, $w \to \infty$.  Then
\eqref{rhoas} follows by Lemma \ref{lemma} and the afore-mentioned properties of $\leb (B_{t/w,r})$.
Theorem  \ref{Xcov} is proved.  \hfill $\Box$

\medskip

In this paper we often use the following lemma which is a version  of Lemma~2 in \cite{kaj2007} or Lemma 2.4 in \cite{bier2010}.

\medskip

\begin{lemma}\label{lemma}
Let $F$ be a probability distribution that has a density function $f$ satisfying \eqref{fbeta}. Set $f_\lambda (r) := \lambda^{1+\alpha} f(\lambda r)$ for $\lambda \ge 1$. Assume that $g$, $g_\lambda$ are measurable functions on $\R_+$ such that $g_\lambda (r) \to g (r)$ as $\lambda \to \infty$ for all $r>0$ and 
such that the 
inequality 
\begin{equation}\label{ineq:bdd_ggc}
| g_\lambda (r) | \ \le \ C (r^{\beta_1} \wedge r^{\beta_2})
\end{equation}
holds for all $r>0$ and some $0 < \beta_1 < \alpha < \beta_2$, where $C$ does not depend on $r, \lambda$. Then
$$
\int_0^\infty g_\lambda (r) f_\lambda (r) \d r \ \to \ c_f \int_0^\infty g (r) r^{-(1+\alpha)} \d r \quad \text{as } \lambda \to \infty.
$$
\end{lemma}

\noi {\it Proof.} Split $\int_0^\infty g_\lambda (r) f_\lambda (r) \d r =
(\int_0^\epsilon  + \int_\epsilon^\infty ) g_\lambda (r) f_\lambda (r) \d r =: I_1(\lambda) + I_2(\lambda) $,
where $\epsilon >0$.
It suffices to prove
\begin{equation} \label{ii}
\lim_{\lambda \to \infty} I_2(\lambda) \ = \ c_f \int_\epsilon^\infty g(r) r^{-(1+\alpha)} \d r
\qquad \text{and} \qquad \lim_{\epsilon \to 0} \limsup_{\lambda \to \infty} I_1(\lambda) \ = \ 0.
\end{equation}
The first relation in \eqref{ii} follows by the dominated convergence theorem, using \eqref{ineq:bdd_ggc} and the bound
$f_\lambda (r) \le C r^{-(1+\alpha)} $ which holds for all $r > \rho/\lambda$ and a sufficiently large $\rho>0$ by
virtue of \eqref{fbeta}. The second relation in \eqref{ii} follows from
$ |I_1(\lambda)|
\le C \int_0^\epsilon r^{\beta_2} f_\lambda (r) \d r = C \lambda^{\alpha-\beta_2} \int_0^{\lambda \epsilon} x^{\beta_2} f(x) \d x
\le  C \lambda^{\alpha-\beta_2}  +   C \lambda^{\alpha-\beta_2}  \int_1^{\lambda \epsilon} x^{\beta_2-(1+\alpha)}
\d x   \le C (\lambda^{\alpha-\beta_2} + \epsilon^{\beta_2 - \alpha}). $ \hfill $\Box$

\subsection{Proofs of Section~4}


\noi {\it Proof of Theorem~\ref{thmslow}.} We have
\begin{eqnarray}
W_{T,\gamma, \beta}(\theta)&:=& \log \E \exp \big\{\i \theta b_T^{-1} \big( A_{M,K}(T x) - \E A_{M,K}(T x) \big) \big\} \nn \\
&=& T^\gamma\int_{\R \times \R_+} \Psi \Big(\theta T^{-\cal H} (r^{1-p} \wedge T^\beta) \int_0^{T x}
\1 ( u < t < u + r^p ) \d t \Big) \d u f(r) \d r, \label{WT}
\end{eqnarray}
where $\Psi(z) = \e^{\i z} - 1 - {\i z}, z \in \R$ as in Sec.~5.1.

\smallskip

\noi (i) Let $0< p < 1,
\delta_1 := {\beta - (1+\gamma)(1-p)/\alpha} > 0, \delta_2 := 1 - (1+\gamma)p/\alpha = (\gamma_+ - \gamma) p / \alpha > 0$.
Using the change of  variables
$(t-u)/r^p \to t$, $u \to T u$, $r \to T^{(1+\gamma)/\alpha} r$ in \eqref{WT},
we obtain
\begin{equation} \label{WTnew}
W_{T,\gamma, \beta}(\theta)\ = \ \int_0^\infty g_T (r) f_T (r) \d r,
\end{equation}
where $f_T (r) := T^{(1+\alpha)(1+\gamma)/\alpha} f (T^{(1+\gamma)/\alpha} r)$ and
$$
g_T (r) \ := \ \int_\R \Psi \big( \theta (r^{1-p} \wedge T^{\delta_1}) r^p h_T (u,r) ) \big) \d u
$$
and where $h_T (u,r) := \int_0^1 \1 ( 0 < u + T^{-\delta_2} r^p t < x ) \d t \to \1(0 < u < x)$ for fixed $(u,r) \in \R \times \R_+$, $u \not \in \{ 0, x \}$.
Hence $g_T (r) \to  g(r) := x \Psi (\theta r) $ follows by the dominated convergence theorem.
The bound $|g_T (r)| \le C \min(r, r^2)$ follows from \eqref{ineq:psi} and $\int_\R h_T (u, r) \d u = x$ with $h_T (u,r) \le 1$.
Finally, by Lemma \ref{lemma}, $W_{T,\gamma, \beta}(\theta)
\to x c_f \int_0^\infty \Psi (\theta r) r^{-(1+\alpha)} \d r = \log \E \exp \{ \i \theta L_\alpha (x,1) \}$, proving part (i) for
$0< p < 1 $. The case $p=1$ follows similarly.

\smallskip

\noi (ii) 
Using the same change of variables as in part (i) we rewrite $ W_{T,\gamma, \beta}(\theta) $ as in \eqref{WTnew},
where
\begin{eqnarray*}
g_T (r) \ := \ \int_{\R} \Psi \big( \theta ( ( T^{-\delta_1} r^{1-p} ) \wedge 1 ) r^p h_T (u,r) \big) \d u,
\end{eqnarray*}
where $\delta_1, f_T (r), h_T (u,r)$ are the same as in    \eqref{WTnew} except that now $\delta_1 <0$. 
Next,  $g_T (r) \to x \Psi ( \theta r^p )$ by the dominated convergence theorem 
while $|g_T (r) | \le C \min (r^p, r^{2p})$ follows by \eqref{ineq:psi} and $\int_\R \min (h_T (u,r), h^2_T (u,r) ) \d u \le C$.
Then  $ W_{T,\gamma, \beta}(\theta) \to  W_{\gamma, \beta}(\theta)
:= x c_f \int_0^\infty \Psi (\theta r^p) r^{-(1+\alpha)} \d r$ follows by Lemma \ref{lemma}. To finish the proof of part (ii) it suffices  
to check that
\begin{eqnarray}\label{def:chf}
 W_{\gamma, \beta}(\theta) \ = \ - x \frac{c_f \Gamma (2-\alpha/p)}{\alpha (1- \alpha/p)} \cos \Big( \frac{\pi \alpha}{2p} \Big) |\theta|^{\alpha/p} 
 \Big( 1 - \i \operatorname{sgn} (\theta) \tan \Big( \frac{\pi \alpha}{2p} \Big) \Big) \ =: \ \log \E \exp \{ \i \theta L_{\alpha/p} (x) \}.
\end{eqnarray}

\smallskip

\noi (iii) Denote $\delta_1 := 1+ \gamma - \alpha \beta/(1-p) > 0, \, \delta_2 := 1- p \beta /(1-p) > 0$. Then by
change of variables: $(t-u)/r^p \to t$, $u \to T u$, $r \to T^{\beta/(1-p)} r$ we rewrite $ W_{T,\gamma, \beta}(\theta) $ as in \eqref{WTnew},
where $f_T (r) := T^{(1+\alpha) \beta/(1-p)} f(T^{\beta/(1-p)} r)$ and
\begin{eqnarray*}
g_T (r) \ := \ \int_\R T^{\delta_1} \Psi \big(\theta T^{-\delta_1/2} (r^{1-p} \wedge 1) r^p h_T (u,r) \big) \d u
\end{eqnarray*}
with $h_T (u,r) := \int_0^1 \1(0 < u +T^{-\delta_2} r^p t < x) \d t \to \1 (0< u < x).$
Then $g_T (r) \to - (\theta^2/2) ( r^{1-p} \wedge 1 )^2 r^{2p} x$ by the
dominated convergence theorem using the bounds
$|\Psi (z)| \le z^2 /2$, $z \in \R$ and $h_T (u,r) \le \1 ( - r^p <  u <  x)$. Moreover,
$|g_T (r)| \le C \min (r^{2p}, r^2)$ holds in view of $\int_\R h^2_T (u,r) \d u \le C$. Using Lemma \ref{lemma} we get
$W_{T,\gamma, \beta}(\theta) \to - (\theta^2/2) x c_f \int_0^\infty (r^{1-p} \wedge 1)^2 r^{2p-(1+\alpha)} \d r = - (\theta^2/2) \sigma^2_1 x,$
where
\begin{equation}\label{const:sigma1}
\sigma^2_1 \ := \ \frac{2 c_f (1-p)}{(2-\alpha)(\alpha-2p)} < \infty
\end{equation}
since $\max(1,2p) < \alpha < 2$. This proves part (iii).

\smallskip

\noi (iv) By the same change of variables as in part (iii), we rewrite $W_{T,\gamma,\beta} (\theta)$ as in \eqref{WTnew}, where
$$
g_T (r) \ := \ \int_{\R} T^{\delta_1} \Psi \big( \theta T^{-\delta_1/2} (\log T)^{-1/2} (r^{1-p} \wedge 1) r^p h_T (u,r) \big) \d u
$$
and $f_T (r)$ and $\delta_1$, $\delta_2 > 0$ and $h_T (u,r) : = \int_0^1 \1 ( 0 < u+ T^{-\delta_2} r^p t < x ) \d t \to \1 ( 0 < u < x )$ are the same as in (iii). We split $W_{T,\gamma,\beta}(\theta) = W_{T,\gamma,\beta}^- (\theta) + W_{T,\gamma,\beta}^0 (\theta) + W_{T,\gamma,\beta}^+ (\theta)$ and next prove that $W_{T,\gamma,\beta}^- (\theta) := \int_0^1 g_T (r) f_T (r) \d r$ and $W_{T,\gamma,\beta}^+ (\theta) := \int_{T^{\delta_1 /2p}}^\infty g_T (r) f_T (r) \d r$ are the remainder terms, whereas 
\begin{eqnarray*}
W_{T,\gamma,\beta}^0 (\theta) \ := \ \int_1^{T^{\delta_1 /2p}} g_T (r) f_T (r) \d r &\sim& - \frac{\theta^2}{2} \frac{x c_f}{\log T} \int_{1}^{T^{\delta_1/2p}} r^{2p-(1+2p)} \d r\\
&=& - \frac{\theta^2}{2} \widehat{\sigma}^2_1 x \ =: \ W_{\gamma,\beta} (\theta),
\end{eqnarray*}
where 
\begin{equation}\label{const:sigma1_hat}
\widehat \sigma_1^2 \ := \ c_f \frac{\delta_1}{2p} \ = \ \frac{c_f}{2p (1-p)} ( (1+\gamma)(1-p) - 2p \beta ).
\end{equation}
By \eqref{fbeta}, there exists $\rho > 0$ such that $f_T (r) \le C r^{-(1+2p)}$ for all $r > \rho/T^{\beta/(1-p)}$. Using this bound along with $\int_{\R} h_T (u,r) \d u = x$, $h_T (u,r) \le 1$ and \eqref{ineq:psi}, we get
\begin{eqnarray*}
| W_{T, \gamma, \beta}^- (\theta) | &\le& \frac{C}{\log T} \int_0^1 r^2 f_T (r) \d r \ = \ O ( (\log T)^{-1} ) \ = \ o(1),\\
\qquad
| W_{T, \gamma, \beta}^+ (\theta) | &\le& C \frac{T^{\delta_1/2}}{(\log T)^{1/2}} \int_{T^{\delta_1/2p}}^\infty r^{p-(1+2p)} \d r \ = \ O ( (\log T)^{-1/2} ) \ = \ o(1).
\end{eqnarray*}
We now consider the main term $W_{T, \gamma, \beta}^0 (\theta)$. Let $\widetilde{W}_{T, \gamma, \beta} (\theta) : = - \frac{\theta^2}{2 \log T} \int_1^{T^{\delta_1/2p}} r^{2p} f_T (r) \d r \int_{\R} h_T^2 (u,r) \d u$. Then, by $|\Psi (z) + z^2/2| \le |z|^3/6$, $z \in \R,$ it follows that 
\begin{eqnarray*}
| W_{T,\gamma,\beta}^0 (\theta) - \widetilde{W}_{T, \gamma,\beta} (\theta) | &\le& \frac{C}{(\log T)^{3/2} T^{\delta_1/2}}\int_1^{T^{\delta_1/2p}} r^{3p} f_T (r) \d r \int_{\R} h_T^3 (u,r) \d u\\
&\le& \frac{C}{(\log T)^{3/2} T^{\delta_1/2}} \int_1^{T^{\delta_1/2p}} r^{p-1} \d r \ = \ O ( (\log T)^{-3/2} ) \ = \ o (1).
\end{eqnarray*}
Finally, we estimate $|\widetilde{W}_{T, \gamma,\beta} (\theta) - W_{\gamma,\beta} (\theta)| \le C (J'_T + J''_T)$, where
\begin{eqnarray*}
J'_T &:=& \frac{1}{\log T} \int_1^{T^{\delta_1/2p}} r^{2p} f_T (r) \d r \int_{\R} | h_T^2 (u,r) - \1 (0 < u < x) | \d u,\\
J''_T &:=& \frac{1}{\log T} \int_1^{T^{\delta_1/2p}} r^{2p} | f_T (r) - c_f r^{-(1+2p)} | \d r.
\end{eqnarray*}
Using 
$$
\int_{\R} | h_T^2 (u,r) - \1 (0 < u < x) | \d u \ \le \ 2 \int_0^1 \d t \int_\R |\1(0 < u + T^{-\delta_2} r^p t < x) - \1 (0 < u < x) | \d u \ \le \ C r^p T^{-\delta_2},
$$
we obtain $J'_T \le C (\log T)^{-1} T^{-\delta_2} \int_1^{T^{\delta_1/2p}} r^{p-1} \d r = o (1)$, since $\delta_1/2 \le \delta_2$ for $\gamma \le \gamma_+$. Then $J''_T = o(1)$ follows from \eqref{fbeta}, since $|f_T (r) - c_f r^{-(1+2p)}| \le \epsilon c_f r^{-(1+2p)}$ for all $r > \rho/ T^{\beta/(1-p)}$ and some $\rho >0$ if given any $\epsilon > 0$. This completes the proof of $W_{T, \gamma, \beta} (\theta) \to - (\theta^2/2) \widehat{\sigma}_1^2 x = \log \E \exp \{ \i \theta \widehat{\sigma}_1 B (x) \}$ as $T \to \infty$ for any $\theta \in \R$. 

\smallskip

\noi (v) After the same change of variables as in part (iii) we get $W_{T, \gamma, \beta} (\theta)$ in \eqref{WTnew}, where 
$$
g_T (r) \ := \ \int_\R \Psi \big( \theta (r^{1-p} \wedge 1) r^p h_T (u,r) \big) \d u
$$
with the same $f_T (r)$ and $h_T (u,t) \to \1(0 < u < x)$ as in (iii). By dominated convergence theorem, $g_T (r) \to x \Psi (\theta (r^{1-p} \wedge 1) r^p)$, where we justify its use by \eqref{ineq:psi}, and $h_T (u,r) \le \1 (-r^p < u < x)$. The bound $|g_T (r)| \le C \min(r^p, r^2)$ follows from \eqref{ineq:psi} and $\int_\R h_T (u, r) \d u = x$ with $h_T (u,r) \le 1$.
Finally, by Lemma \ref{lemma},
\begin{eqnarray}\label{def:intlevy}
W_{T, \gamma, \beta} (\theta) \ \to \ x c_f \int_0^\infty \Psi \big( \theta ( r^{1-p} \wedge 1 ) r^p \big) r^{-(1+\alpha)} \d r \ =: \ \log \E \exp \{ \i \theta \widehat L (x) \}. 
\end{eqnarray}
The proof of Theorem \ref{thmslow} is complete.
\hfill $\Box$

\bigskip

\noi {\it Proof of Theorem \ref{thmfast}}. (i) Denote $\delta_1 := 1 + \gamma - \alpha/p = \gamma - \gamma_+ > 0$ and $\delta_2 := (1-p)/p - \beta > 0$.
By changing the variables in \eqref{TH}: $t \to T t$, $u \to T u$, $r \to T^{1/p} r$
we rewrite $ W_{T,\gamma, \beta}(\theta) $ as in \eqref{WTnew}, where
$f_T (r) := T^{(1+\alpha)/p} f (T^{1/p} r)$ and
$$
g_T (r) \ := \ \int_\R T^{\delta_1} \Psi \big( \theta T^{-\delta_1/2} ( (T^{\delta_2} r^{1-p} ) \wedge 1 ) h (u,r) \big) \d u
$$
with $h (u,r) := \int_0^x \1 (u < t < u + r^p  ) \d t$. The dominated convergence $g_T (r) \to  g(r) := - (\theta^2/2) \int_\R h^2 (u,r) \d u$ follows by \eqref{ineq:psi}. The latter combined with $\int_\R h^2 (u,r) \d u \le C \min (1,r^p) \int_\R h (u,r) \d u \le C \min (r^p, r^{2p})$ gives the bound $|g_T (r) | \le C \min (r^p, r^{2p})$. Finally, by Lemma \ref{lemma}, $W_{T,\gamma, \beta}(\theta) \to - (\theta^2/2) \sigma^2_2 x^{2 H}$, where
\begin{eqnarray}\label{const:sigma2}
\sigma^2_2 &:=& c_f \int_{\R \times \R} \Big( \int_0^1 \1 ( u < t < u + r^p) \d  t \Big)^2 \frac{\d u \d r}{r^{1+\alpha}}  \
=\  \frac{2 c_f}{\alpha (2 - \alpha/p)(3-\alpha/p)(\alpha /p - 1)},
\end{eqnarray}
proving part (i).

\smallskip

\noi (ii) The proof is the same as of Theorem \ref{thmslow} (iii). 

\smallskip

\noi (iii) Let $\delta_1 := \gamma - \alpha_+ \beta > 0, \, \delta_2 := \alpha_+ \beta/\gamma_+ -1> 0$. By
change of variables: $t \to T t$, $u \to T^{\beta p/(1-p)} u$, $r \to T^{\beta/(1-p)} r$ we get  \eqref{WTnew} with
$f_T (r) := T^{(1+\alpha)\beta / (1-p)} f (T^{\beta / (1-p)} r)$ and
$$
g_T (r) \ := \ \int_\R T^{\delta_1} \Psi ( \theta T^{-\delta_1/2} ( r^{1-p} \wedge 1 ) h_T (u,r) ) \d u,
$$
with $h_T (u,r) := \int_0^x \1 (0 < (T^{-\delta_2} t - u)/r^p < 1) \d t \to  h(u,r) := x \1 (-r^p < u < 0) $.
Then \eqref{ineq:psi} and $h^2_T (u,r) \le x \1 ( -r^p < u < 1)$ justify the dominated convergence $g_T (r) \to - (\theta^2/2) (r^{1-p} \wedge 1)^2  r^p x^2$. By \eqref{ineq:psi} and $\int_\R h^2_T (u,r) \d u \le C \int_\R h_T (u,r) \d u  \le C r^p$, we have $|g_T (r)| \le C \min (r^p, r^{2-p})$. 
Finally, by Lemma \ref{lemma} 
$W_{T,\gamma, \beta}(\theta)  \to - (\theta^2/2) x^2 c_f \int_0^\infty (r^{1-p} \wedge 1)^2 r^{p-(1+\alpha)} \d r = - (\theta^2/2) x^2 \sigma^2_3$ with
\begin{equation}\label{const:sigma3}
\sigma^2_3 \ := \ \frac{2 c_f (1-p)}{(2 - p - \alpha) (\alpha - p)},
\end{equation}
proving part (iii).

\smallskip

\noi (iv) Denote
$\delta_1 := \beta - \gamma/\alpha_+ >0,
\delta_2 := \gamma/\gamma_+ - 1 >0$.  Using the change of variables: 
$t \to T t$, $u \to T^{\gamma/\gamma_+} u$, $r \to T^{\gamma/\gamma_+ p} r$ we get  \eqref{WTnew} with
$f_T (r) := T^{(1+\alpha)\gamma/\gamma_+ p} f (T^{\gamma/\gamma_+ p} r)$ and
\begin{eqnarray*}
g_T (r) \ := \ \int_\R \Psi ( \theta ( r^{1-p} \wedge T^{\delta_1} ) h_T (u,r) ) \d u, 
\end{eqnarray*}
where $h_T(u,r) := \int_0^x \1 (u < T^{-\delta_2} t < u + r^p) \d t \to h(u,r) := x \1( -r^p < u < 0)$.
Then   $g_T (r) \to g(r) :=  \int_\R \Psi ( \theta x r^{1-p} \1( -r^p < u < 0) ) \d u $ and
$W_{T,\gamma, \beta}(\theta)  \to  c_f \int_0^\infty  g(r) r^{-(1+\alpha)} \d r 
=  \log \E \exp \{\i \theta x L_+(1)\} $ similarly as in the proof of Theorem \ref{thm2} (ii).

\smallskip

\noi (v) Set $\delta_1 := \gamma - \gamma_+ > 0$, $\delta_2 := \beta - (1-p)/p > 0$. After a change of variables: $t \to T t$, $u \to T u$, $r \to T^{1/p} r$, we get \eqref{WTnew} with $f_T (r) := T^{(1+\alpha)/p} f(T^{1/p} r)$ and
$$
g_T (r) \ := \ \int_\R T^{\delta_1} \Psi ( \theta T^{-\delta_1/2} (r^{1-p} \wedge T^{\delta_2}) h (u,r)) \d u ,
$$
where $h(u,r) := \int_0^x \1 (u < t < u + r^p) \d t$. Then $g_T (r) \to g (r) := - (\theta^2/2) \int_{\R} r^{2(1-p)} h^2(u,r) \d u$ and $W_{T, \gamma, \beta} (\theta) \to c_f \int_0^\infty g(r) r^{-(1+\alpha)} \d r = - (\theta^2/2) \sigma^2_+ x^{2 H_+}$ similarly to the proof of Theorem \ref{thm3} (i). 

\smallskip

\noi (iii) We follow the proof of Theorem \ref{thmslow} (iv). By the same change of variables, we rewrite $W_{T,\gamma,\beta} (\theta)$ as in \eqref{WTnew}. We split $W_{T,\gamma,\beta} (\theta) = W_{T,\gamma,\beta}^- (\theta) + W_{T,\gamma,\beta}^0 (\theta) + W_{T,\gamma,\beta}^+ (\theta)$ with the same $W_{T,\gamma,\beta}^\pm (\theta)$ being the remainder terms. Note that now $\delta_2 < \delta_1/2$, since $\gamma > \gamma_+$. Next, we split $W_{T,\gamma,\beta}^0 (\theta) = W'_{T,\gamma,\beta} (\theta) + W''_{T,\gamma,\beta} (\theta)$, where 
$$
W'_{T,\gamma,\beta} (\theta) \ := \ \int _1^{T^{\delta_2/p}} g_T (r) f_T (r) \d r, \qquad W''_{T,\gamma,\beta} (\theta) \ := \ \int_{T^{\delta_2/p}}^{T^{\delta_1/2p}} g_T (r) f_T (r) \d r.
$$
Analogously to the proof of Theorem \ref{thmslow} (iv), we show the convergence $W'_{T,\gamma,\beta} (\theta) \to -(\theta^2/2) \widehat{\sigma}^2_2 x$, where 
\begin{equation}\label{const:sigma_hat_2}
\widehat \sigma^2_2 \ := \ c_f \frac{\delta_2}{p} \ = \ c_f \Big( \frac{1}{p} - \frac{\beta}{1-p} \Big).
\end{equation}
Using \eqref{ineq:psi} and $\int_{\R} h_T (r,u) \d u = x$ with $h_T (r,u) \le x (T^{\delta_2} / r^p)$, we get 
$$
|  W''_{T,\gamma,\beta} (\theta) | \ \le \ \frac{C}{\log T} \int_{T^{\delta_2/p}}^{T^{\delta_1/2p}} \frac{\d r}{r} \int_{\R}  h_T^2 (r, u) \d u \ \le \  \frac{C  T^{\delta_2}}{\log T} \int_{T^{\delta_2 / p}}^{T^{\delta_1/2p}} \frac{\d r}{r^{1+p}} \ = \ O ( (\log T)^{-1} ) \ = \ o (1),
$$
which completes the proof of $W_{T,\gamma,\beta} (\theta) \to -(\theta^2/2) \widehat{\sigma}^2_2 x = \log \E \exp \{ \i \theta \widehat{\sigma}_2 B (x) \}$ as $T \to \infty$ for any $\theta \in \R$. 

\smallskip

\noi (vii) By the same change of variables as in part (i), we rewrite $W_{T,\gamma,\beta} (\theta)$ as in \eqref{WTnew}, where 
$$
g_T (r) \ := \ \int_{\R} T^{\delta_1} \Psi \big( \theta T^{-\delta_1/2} (r^{1-p} \wedge 1) h(u,r) \big) \d u
$$
and where $\delta_1$, $h(u,r)$, $f_T(r)$ are the same as in (i). Then $g_T (r) \to - (\theta^2 / 2) \int_{\R} h^2  (u,r) \d u$ along with $\int_R h^2(u,r) \d u \le C \min (r^p, r^{2p})$ and \eqref{ineq:psi} imply $W_{T,\gamma,\beta} (\theta) \to - (\theta^2/2) c_f \int_0^\infty \int_{\R} (r^{1-p} \wedge 1)^2 h^2(u,r) r^{-(1+\alpha)} \d r \d u =: \log \E \exp \{ \i \theta \widehat{Z} (x) \}$ as $T \to \infty$ for any $\theta \in \R$, by Lemma \ref{lemma}.

\smallskip

\noi (viii) By the same change of variables as in part (iii) we obtain $W_{T,\gamma,\beta} (\theta)$ as in \eqref{WTnew}, where $g_T (r) := \int_{\R} \Psi ( \theta (r^{1-p} \wedge 1) h_T (u,r) ) \d u$ with $f_T (r)$, $\delta_2 = \gamma / \gamma_+ - 1 > 0$ and $h_T (u,r) := \int_0^x \1 ( u < T^{-\delta_2} t < u + r^p ) \d t \to x \1 ( -r^p < u < 0 )$ the same as in (iii). 
Using $\int_{\R} h_T (u,r) \d u = x r^p$ and $h_T (u,r)\le x$ yields $|g_T (r) | \le C \min (r^{p}, r^{2-p})$ from \eqref{ineq:psi}. Hence, by Lemma \ref{lemma}, it follows that
\begin{eqnarray}\label{def:rv_hatZ}
W_{T,\gamma,\beta} (\theta) \ \to \ c_f \int_0^\infty \Psi ( \theta x (r^{1-p} \wedge 1) ) r^{p-(1+\alpha)}\d r \ =: \ \log \E \exp \{ \i \theta x \widehat Z \}.
\end{eqnarray}

\smallskip

\noi (ix) By the same change of variables as in the proof of part (iv), we rewrite $W_{T,\gamma,\beta} (\theta)$ as in \eqref{WTnew}, where 
$$
g_T (r) \ := \ \int_\R \Psi \big( \theta (\log T)^{-1/2} (r^{1-p} \wedge T^{\delta_1}) h_T (u,r) \big) \d u
$$
with $\delta_1$, $\delta_2 := \gamma/\gamma_+ - 1 > 0 $ and $h_T (u,r) := \int_0^x \1 ( u < T^{-\delta_2} t < u + r^p ) \d t \to x \1 ( -r^p < u < 0 ) =: h (u,r)$ and $f_T (r)$ being the same as in (iv). We split $W_{T,\gamma,\beta} (\theta) = W^-_{T,\gamma,\beta} (\theta) + W^0_{T,\gamma,\beta} (\theta) + W^+_{T,\gamma,\beta} (\theta)$ and next prove that $W^-_{T,\gamma,\beta} (\theta) := \int_0^{T^{-\delta_2/p}} g_T (r) f_T (r)$ and $W^+_{T,\gamma,\beta} (\theta) := \int_1^\infty g_T (r) f_T (r) \d r$ are the remainder terms, whereas 
\begin{eqnarray*}
W^0_{T,\gamma,\beta} (\theta) \ := \ \int_{T^{-\delta_2/p}}^1 g_T (r) f_T (r) \d r &\sim& - \frac{\theta^2}{2} \frac{c_f}{\log T} \int_{T^{-\delta_2/p}}^1 \frac{\d r}{r^{1+p}} \int_\R h^2 (u,r) \d u\\
&=& - \frac{\theta^2}{2} \widetilde{\sigma}^2_+ x^2 \ =: \ W_{\gamma,\beta} (\theta),
\end{eqnarray*}
where the constant $\widetilde{\sigma}^2_+$ is given in \eqref{wtisigma}.
Using $\int_\R h_T (u,r) \d u = x r^p$ and $h_T (u,r) \le x \wedge (T^{\delta_2} r^p)$ along with \eqref{ineq:psi}, we show that
\begin{eqnarray*}
| W^-_{T,\gamma,\beta} (\theta) | &\le& \frac{C T^{\delta_2}}{\log T} \int_0^{T^{-\delta_2/p}} r^2 f_T (r) \ = \ \frac{C}{T \log T} \int_0^{T^{1/p}} r^2 f(r) \d r \ = \ O ( (\log T)^{-1} ) \ = \ o(1),\\
| W^+_{T,\gamma,\beta} (\theta) | &\le& \frac{C}{(\log T)^{1/2}} \int_1^\infty r f_T (r) \d r \ = \ O ( (\log T)^{-1/2} ) \ = \ o (1).
\end{eqnarray*}
To deal with the main term $W^0_{T,\gamma,\beta} (\theta)$, set $\widetilde{W}_{T,\gamma,\beta} (\theta) := - \frac{\theta^2}{2 \log T} \int_{T^{-\delta_2/p}}^1 r^{2(1-p)} f_T (r) \d r \int_\R h_T^2 (u,r) \d u$. From $| \Psi (z) + z^2/2 | \le |z|^3/6$, we obtain
\begin{eqnarray*}
| W_{T,\gamma,\beta} (\theta) - \widetilde{W}_{T,\gamma,\beta} (\theta) | &\le& \frac{C}{(\log T)^{3/2}} \int_{T^{-\delta_2/p}}^1 r^{3(1-p)} f_T (r) \d r \int_\R h_T^3 (u,r) \d u\\
 &\le& \frac{C}{(\log T)^{3/2}} \int_{T^{-\delta_2/p}}^1 r^{3-2p} f_T (r) \d r \ = \ O ((\log T)^{-3/2}) \ = \ o(1).
\end{eqnarray*}
Finally, we consider $| \widetilde{W}_{T,\gamma,\beta} (\theta) - W_{\gamma,\beta} (\theta) | \le C (J'_T + J''_T)$, where
\begin{eqnarray*}
J'_T &:=& \frac{1}{\log T} \int_{T^{-\delta_2/p}}^1 r^{2(1-p)} f_T (r) \d r \int_\R | h_T^2 (u,r) - h^2 (u,r)  | \d u,\\
J''_T &:=& \frac{1}{\log T} \int_{T^{-\delta_2/p}}^1 r^{2-p} | f_T (r) - c_f r^{p-3} | \d r.
\end{eqnarray*}
Using 
$$
\int_\R | h^2_T (u,r) - h^2 (u,r) | \d u \ \le \ C \int_0^x \d t \int_\R | \1 (u < T^{-\delta_2} t < u + r^p ) - \1 (-r^p < u < 0)) | \d u \ \le \ C T^{-\delta_2}
$$
we obtain $J'_T \le C (\log T)^{-1} T^{-\delta_2} \int_{T^{-\delta_2/p}}^1 r^{-(1+p)} \d r = O( (\log T)^{-1} ) = o (1)$.  
Then $J''_T = o (1)$ follows from \eqref{fbeta}, since $| f_T (r) - c_f r^{p-3} | \le \epsilon c_f r^{p-3}$ for all $r > \rho/T^{\gamma/2(1-p)}$ and some $\rho > 0$ if given any $\epsilon > 0$. This finishes the proof of $W_{T,\gamma,\beta} (\theta) \to - (\theta^2/2) \widetilde{\sigma}^2_+ x^2 = \log \E \exp \{ \i \theta \widetilde{\sigma}^2_+ B_{1,1/2} (x,1) \} $ as $T \to \infty$ for any $\theta \in \R$. 

\noi The proof of Theorem \ref{thmfast} is complete.
\hfill $\Box$

\bigskip

\noi {\it Proof of Theorem \ref{thminter}.} (i) By the same change of variables as in Theorem \ref{thmfast} (i), we rewrite $W_{T,\gamma,\beta} (\theta)$ as in \eqref{WTnew}, where 
$$
g_T (r) \ := \ \int_{\R} \Psi \big( \theta  ( ( T^{\delta_2} r^{1-p} ) \wedge 1 ) h(u,r) \big) \d u \ \to \ \int_{\R} \Psi ( \theta h(u,r) ) \d u \ =: \ g (r),
$$
since $\delta_2 := (1-p)/p - \beta = \gamma_+/\alpha_+ - \beta > 0$ with $h(u,r)$, $f_T (r)$ being the same as in Theorem~\ref{thmfast}~(i).
Using \eqref{ineq:psi} along with $\int_\R h (u,r) \d u = x r^p$ and $h (u,r) \le r^p$, we get $|g_T (r) | \le C \min( r^p, r^{2p} )$. Hence $W_{T,\gamma,\beta} (\theta) \to c_f \int_0^\infty \int_\R \Psi ( \theta h (u,r) ) r^{-(1+\alpha)} \d r \d u =: \log \E \exp \{ \i \theta I (x) \}$ by Lemma \ref{lemma}.

\smallskip

\noi (ii), (iii) The proof is the same as that of Theorem \ref{thmslow} (iii), (iv) respectively.

\smallskip

\noi (iv) By the same change of variables as in Theorem \ref{thmfast} (i), we rewrite $W_{T,\gamma,\beta} (\theta)$ as in \eqref{WTnew}, where 
$
g (r) := \int_{\R} \Psi ( \theta (r^{1-p} \wedge 1) h(u,r) ) \d u
$
with $h(u,r)$, $f_T(r)$ being the same as in Theorem \ref{thmfast} (i).
Then $|g (r)| \le C \min (r^p, r^2)$ follows from \eqref{ineq:psi}. By Lemma \ref{lemma}, we get $W_{T,\gamma,\beta} (\theta) \to c_f \int_0^\infty g(r) r^{-(1+\alpha)} \d r =: \log \E \exp \{ \i \theta \widehat I (x) \}$.

\smallskip

\noi (v) By the same change of variables as in Theorem \ref{thmfast} (v), we rewrite $W_{T,\gamma,\beta} (\theta)$  as in \eqref{WTnew}, where $f_T(r)$, $g_T(r)$ are the same as in Theorem \ref{thmfast} (v) except for $\delta_1 = 0$. Then $g_T (r) \to g(r) := \int_\R \Psi (\theta r^{1-p} h(u,r)) \d u$ and $|g_T(r)| \le C \min (r,r^2)$ from \eqref{ineq:psi} lead to $W_{T,\gamma,\beta} (\theta) \to c_f \int_0^\infty g(r) r^{-(1+\alpha)} \d r = \log \E \exp \{ \i \theta I_+ (x,1) \}$ by Lemma~\ref{lemma}, similarly to the proof of Theorem \ref{thm5}. 

\noi The proof of Theorem \ref{thminter} is complete. 
\hfill $\Box$

\section*{Acknowledgement}

This research was supported by a grant (no. MIP-063/2013) from the Research Council of Lithuania.

\end{document}